\magnification=1200
\hsize=125mm       \vsize=187mm
\hoffset=4mm       \voffset=10mm
\pretolerance=500  \tolerance =1000  \brokenpenalty=5000

\parindent=10mm       \hfuzz=1pt      \parskip=4mm

\font\pcap=cmcsc10
\font\pcapg=cmcsc10 scaled 1200
\font\fourteenbf=cmbx10 at 14pt
\font\sixteenrm=cmr10 at 16pt

\catcode`\@=11
\catcode`\;=\active
\def;{\relax\ifhmode\ifdim\lastskip>\z@
\unskip\fi\kern.2em\fi\string;}
\catcode`\:=\active
\def:{%\relax\ifhmode\ifdim\lastskip>\z@
%\unskip\fi\penalty\@M\ \fi
\string:}
\catcode`\!=\active
\def!{\relax\ifhmode\ifdim\lastskip>\z@
\unskip\fi\kern.2em\fi\string!}
\catcode`\?=\active
\def?{\relax\ifhmode\ifdim\lastskip>\z@
\unskip\fi\kern.2em\fi\string?}
\def\^#1{\if#1i{\accent"5E\i}\else{\accent"5E #1}\fi}
\def\"#1{\if#1i{\accent"7F\i}\else{\accent"7F #1}\fi}
\catcode`\@=12
\frenchspacing
\def\og{\leavevmode\raise.3ex\hbox{$\scriptscriptstyle\langle\!\langle$}}
\def\fg{\leavevmode\raise.3ex\hbox{$\scriptscriptstyle\,\rangle\!\rangle$}}

\def\NM{{\rm{I\kern -.2em {\bf N}}}}
\def\CM{{\rm\rlap{\kern 0.3em l}{\bf C}}}
\def\RM{{\rm{I}\kern -.2em{\bf R}}}

\def\HM{{\rm{I\kern -.2em {\bf H}}}}

\long\def\titre#1{{\vskip 1cm\noindent{\fourteenbf #1.}}}
\long\def\theo#1#2{\medskip\begingroup\noindent{\bf#1.\hskip 2mm}%
{\it#2}\endgroup}
\long\def\rmq#1#2{\medskip\noindent{\pcap#1: }{#2}}
\def\ref#1{#1}
\def\inter{\rlap{\raise 8pt\hbox{$\scriptstyle\,\,\,\circ$}}}

\input psfig.sty
\long\def\figure#1#2{\vbox{\noindent\leavevmode\centerline{%
\psfig{file=#1}}\smallskip\centerline{#2}}}
%%%%%%%%%%%%%%%%%%%%%%%%%%%%%%%%%%%%%%%%%%%%%%%%
\null
{\sixteenrm \centerline{Construction of harmonic diffeomorphisms}
\vskip 4mm
\centerline{and minimal graphs}\par}
\vskip 6mm
\centerline{\pcapg By Pascal Collin and Harold Rosenberg }
\centerline{\pcap (version1 / January~2007)}
\vskip10mm
\begingroup
\setbox100=\vbox{\hsize=115mm\noindent{\bf Abstract:} We study complete minimal graphs in $\HM\times\RM$, which take asymptotic boundary values plus and minus infinity on alternating sides of an ideal inscribed polygon $\Gamma$ in $\HM$.  We give necessary and sufficient conditions on the "lenghts" of the sides of the polygon (and all inscribed polygons in $\Gamma$) that ensure the existence of such a graph.  We then apply this to construct entire minimal graphs in $\HM\times\RM$ that are conformally the complex
plane $\CM$. The vertical projection of such a graph yields a harmonic diffeomorphism from $\CM$ onto $\HM$, disproving a conjecture of Rick Schoen.\par
}
\hbox{\centerline{\box100}}
\endgroup

\titre{1.~Introduction}

In 1952, E.~Heinz proved there is no harmonic diffeomorphism from a disk onto the complex plane $\CM$, with the euclidean metric [\ref{He}]. He used this to give another proof of Bernstein's theorem: an entire minimal graph over the euclidean plane is a plane.

Later, R.~Schoen and S.T.~Yau, asked whether Riemannian surfaces which are related by a harmonic diffeomorphism are quasi-conformally related. In that direction, R.~Schoen conjectured there is no harmonic diffeomorphism from $\CM$ onto the hyperbolic plane $\HM$ [\ref{S}], [\ref{S-Y}] and [\ref{M}].

In this paper we will construct harmonic diffeomorphisms from $\CM$ onto $\HM$. We will use entire minimal graphs to construct these examples (E.~Heinz used the non-existence of harmonic diffeomorphisms from $\HM$ onto $\CM$ to prove the non-existence of non-trivial entire minimal euclidean graphs).

Consider the Riemannian product $\HM\times\RM$ and entire minimal graphs $\Sigma$ defined over $\HM$. The vertical projection $\Sigma\longrightarrow\HM$ is a surjective harmonic diffeomorphism, so we will solve the problem by constructing entire minimal graphs $\Sigma$ that are conformally $\CM$ (Theorem~3). Notice that the horizontal projection $\Sigma\longrightarrow\RM$ is a harmonic function on $\Sigma$ when $\Sigma$ is a minimal surface. Hence if this height function is bounded on $\Sigma$, $\Sigma$ is necessarily hyperbolic (conformally the unit disc). So we must look for unbounded minimal graphs.

Here is the idea of the construction. Let $\Gamma$ be an ideal geodesic polygon in $\HM$ with an even number of sides (the vertices of $\Gamma$ are at infinity). We give necessary and sufficient conditions on the geometry of $\Gamma$ which ensure the existence of a minimal graph $u$ over the polygonal domain $D$ bounded by $\Gamma$, which takes the values plus and minus infinity on alternate sides of $\Gamma$ (\S~3). This is a Jenkins-Serrin type theorem at infinity [\ref{J-S}]. We call such graphs $u$ over $D$ ideal Scherk graphs and we show their conformal type is $\CM$ (\S~5).

We attach certain ideal quadrilaterals to all of the sides of $\Gamma$ (outside of $D$) so that the extended polygonal domain $D_1$ admits an ideal Scherk graph. We do this so that $D_1$ depends on a small parameter $\tau$, and a minimal Scherk function $u_1(\tau)$ defined over $D_1$ satisfies the following. Given a fixed compact disk $K_0$ in the domain $D$, $u_1( \tau)$ is as close as we wish to $u$ over $K_0$, for $\tau$ sufficiently small. Also ideal Scherk graphs are conformally $\CM$ so there is a compact disk $K_1$ ( containing $K_0$) in each $D_1$, so that the conformal type of the annulus in the graph of $u_1(\tau)$, over $K_1-K_0$, is greater than one. Now fix $\tau$ and do the same process to enlarge $D_1$, attaching certain ideal quadrilaterals to all of the sides to obtain a domain, admitting an ideal Scherk graph $u_2(\tau')$, which is as close as we want to $u_1(\tau)$ on $K_1$ for $\tau'$ sufficiently small; then the conformal moduli of the annulus in the graph of $u_2(\tau')$ over $K_1-K_0$ remains greater than one. As before take $K_2$ with the same condition on modulus for $u_2(\tau')$ over $K_2-K_1$. The entire minimal graph $\Sigma$ is obtained by continuing this process and choosing a convergent subsequence. The conformal type of $\Sigma$ is $\CM$ because we write $\HM$ as the union of an increasing sequence of compact disks $K_n$, and each annulus in $\Sigma$ over each $K_{n+1}-K_n$ has conformal modulus at least one.

\titre{2.~Structure of the divergence set}

First we will state some properties of solutions established in [\ref{J-S}] and [\ref{N-R}]. By solution in $D$ we mean a solution of the minimal surface equation in a domain $D$ of the hyperbolic plane. We make no assumptions here on $\partial D$. Starting with curvature estimates for stable minimal surfaces, we give new proofs of the results of this section in an appendix.

\theo{Compactness Theorem}{Let $\{u_n\}$ be a uniformly bounded sequence of solutions in $D$. Then a subsequence converges uniformly on compact subsets of $D$, to a solution in $D$.}

\theo{Monotone Convergence Theorem}{Let $\{u_n\}$ be a monotone sequence of solutions in $D$. If the sequence $\{\vert u_n\vert\}$ is bounded at one point of $D$, then there is a non-empty open set $U\subset D$ (the convergence set) such that $\{u_n\}$ converges to a solution in $U$. The convergence is uniform on compact subsets of $U$ and the divergence is uniform on compact subsets of $D-U=V$. $V$ is called the divergence set.}

Now assume $\partial D$ is an ideal polygon with a finite number of vertices at $\partial_\infty\HM$, composed of geodesic arcs $A_1,\ldots,A_k$, $B_1,\ldots,B_{k'}$ joining the vertices, together with convex arcs (convex towards $D$), $C_1,\ldots,C_{k''}$. We assume $D$ is simply connected, $\partial D$, together with the vertices, homeomorphic to $S^1$, and no two $A$'s (or $B$'s) have a vertex in common.

\theo{Divergence Structure Theorem}{Let $\{u_n\}$ be a monotone sequence of solutions in $D$, each $u_n$ continuous on $\overline D$. If the divergence set $V\neq\emptyset$, then $\rm{int}\,V\neq\emptyset$, and $\partial V$ is composed of ideal geodesics among the $A_i$ and $B_j$, and interior ideal geodesics $C\subset D$, joining two vertices of $\partial D$. No two interior geodesics $C_1,C_2$ of $\partial V$ go to the same vertex at infinity.}

\rmq{Remark~1}{With the exception of the last sentence in the above theorem the proofs are in the papers cited at the beginning of this section. We will prove the last statement after stating the flux relations.}

Now let $\{u_n\}$ be a sequence defined in $D$ satisfying the hypothesis of the Divergence Structure Theorem. For each $n$, let $X_n={\nabla u_n\over W_n}$ be the vector field on $D$, $W_n^2=1+\vert\nabla u_n\vert^2$. For ${\cal W}\subset D$, and $\alpha$ a boundary arc of ${\cal W}$, we define the flux of $u_n$ across $\alpha$ to be $F_n(\alpha)=\int_\alpha\langle X_n,\nu\rangle\,ds$; here $\alpha$ is oriented as the boundary of ${\cal W}$ and $\nu$ is the outer conormal to ${\cal W}$ along $\alpha$. More generaly, for any solution $u$ in $D$ and an oriented arc $\alpha$, we write 
$F_u(\alpha)$ the flux of the associate field $X={\nabla u\over W}$, $W=(1+\vert\nabla u\vert^2)^{1/2}$.

\theo{Flux Theorem}{Let ${\cal W}$ be a domain in $D$. Then 

\item{\it i)}if $\partial{\cal W}$ is a compact cycle, then $F_n(\partial{\cal W})=0$;
\item{\it ii)}if ${\cal W}\subset U$ (the convergence set) and $\alpha$ is a compact arc of 
$\partial{\cal W}$ on which the $u_n$ diverge to $+\infty$, then $\alpha$ is a geodesic and
$$\lim_{n\to\infty} F_n(\alpha)=\vert\alpha\vert.$$
If the $u_n$ diverge to $-\infty$ on $\alpha$, then $\alpha$ is a geodesic and
$$\lim_{n\to\infty} F_n(\alpha)=-\vert\alpha\vert.$$
\item{\it iii)}if ${\cal W}\subset V$, and the $u_n$ remain uniformly bounded on $\alpha$, then
$$\lim_{n\to\infty} F_n(\alpha)=-\vert\alpha\vert,\hbox{ if }u_n\to+\infty\hbox{ on }V,$$
and
$$\lim_{n\to\infty} F_n(\alpha)=\vert\alpha\vert,\hbox{ if }u_n\to-\infty\hbox{ on }V.$$}

\rmq{Remark~2}{In {\it ii} of the Flux Theorem, the fact that $\alpha$ (contained in $\partial{\cal W}$) is a geodesic when the $u_n$ diverge on $\alpha$, is called the Straight Line Lemma.}

Now we can prove the last statement of the Divergence Structure Theorem.

Suppose on the contrary, that $\partial V$ has two interior arcs $C_1$, $C_2$ going to the same vertex at infinity of $D$. Let ${\cal W}\subset D$ be a bounded domain, with $\partial{\cal W}$ a simple closed curve composed of two horocycle arcs $\gamma_1$, $\gamma_2$ and two arcs $\alpha_1$, $\alpha_2$ contained in $C_1$, $C_2$ respectively.

By part {\it i} of the Flux Theorem, $F_n(\partial{\cal W})=0$, i.e.,
$$F_n(\alpha_1)+F_n(\gamma_1)+F_n(\alpha_2)+F_n(\gamma_2)=0.$$
First assume the $u_n$ diverge to $+\infty$ on ${\cal W}$; i.e. ${\cal W}\subset V$, cf. figure~1. Since $U$ is on the other side from ${\cal W}$ along $\alpha_1$ and $\alpha_2$, we have
$$\lim_{n\to\infty}F_n(\alpha_1)=-\vert\alpha_1\vert=\lim_{n\to\infty} F_n(\alpha_2),$$
by part {\it iii} of the Flux Theorem. But $F_n(\gamma_i)\leq\vert\gamma_i\vert$ for $i=1,2$, and the length of $\gamma_2$ can be chosen arbitrarily small, the length of the $\alpha_i$ arbitrarily large, which contradicts the flux equality along $\partial{\cal W}$.

If the $u_n$ diverge to $-\infty$ in ${\cal W}$, one obtains a similar contradiction.

\midinsert
\figure{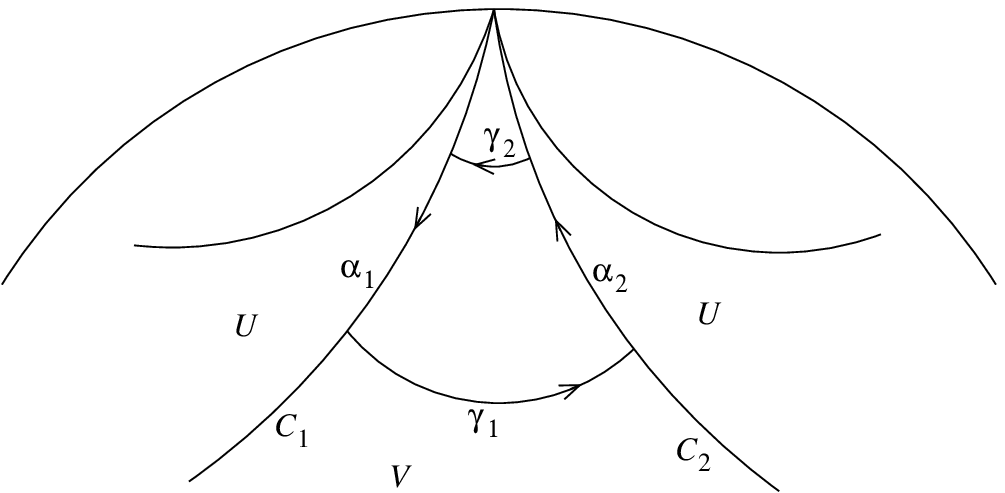}{\rlap{\smash{\raise25mm\hbox{\hskip3mm${\cal W}$}}}Figure~1}
\endinsert

When ${\cal W}$ is not contained in $V$, then ${\cal W}\subset U$ and the same flux equation gives a contradiction; the only change is the sign of 
$\lim_{n\to\infty}F_n(\alpha_1)=\lim_{n\to\infty} F_n(\alpha_2)$.

When we establish existence theorems for unbounded boundary data, we need to know solutions take on the boundary values prescribed. This is guaranteed in our situation by the following result.

\theo{Boundary Values Lemma}{Let $D$ be a domain and let $C$ be a compact convex arc in $\partial D$. Suppose $\{u_n\}$ is a sequence of solutions in $D$ that converges uniformly on compact subsets of $D$ to a solution $u$ in $D$. Assume each $u_n$ is continuous in $D\cup C$ and that the boundary values of $u_n$ on $C$ converge uniformly to a function $f$ on $C$. Then $u$ is continuous in $D\cup C$ and $u$ equals $f$ on $C$.}

\noindent{\bf Proof:} One needs to show that for each $p\in C$, the sequence $\{u_n\}$ is uniformly bounded in a neighborhood of $p$ in $D\cup C$; then standard local barriers show $U$ takes on the correct boundary values.

When $C$ is strictly convex in a neighborhood of $p$ then one places a Scherk surface defined over a geodesic triangle over the graph of $u$; cf. figure~2.

\midinsert
\figure{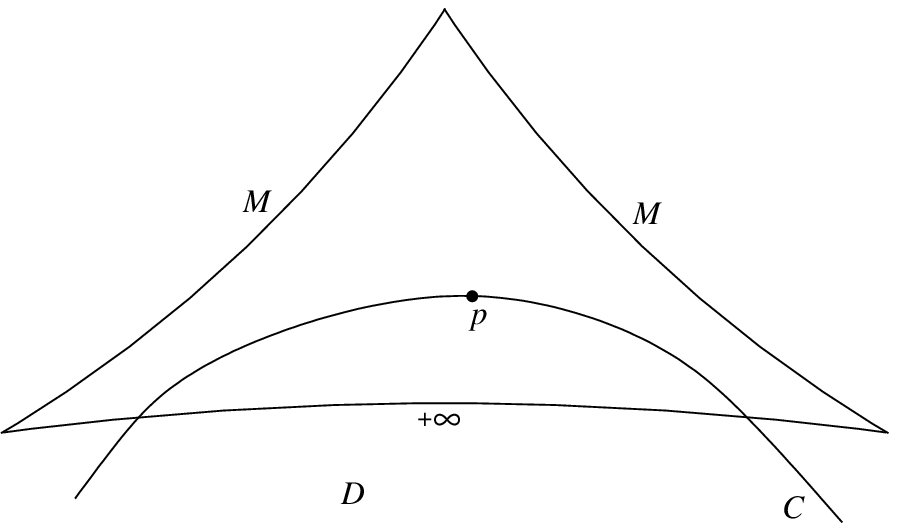}{Figure~2}
\endinsert

Here $M$ is the maximum of $f$ on $C$. Since the boundary values of this Scherk surface are above those of $u_n$, the Scherk surface is above the graph of $u_n$ on the region in the triangle over $\overline D$. Thus the Scherk surface is above $u$, and this uniformly bounds $u$ in a smaller compact neighborhood of $p$ inside the triangle.

When the arc $C$ contains a geodesic segment in a neighborhood of $p$, one uses an analogous Scherk barrier.

Consider a quadrilateral ${\cal P}$ composed of two (``vertical'') geodesics $A_1$, $A_2$, and two geodesic curves $C_1$, $C_2$ (``horizontal''); cf. figure~3.

\midinsert
\figure{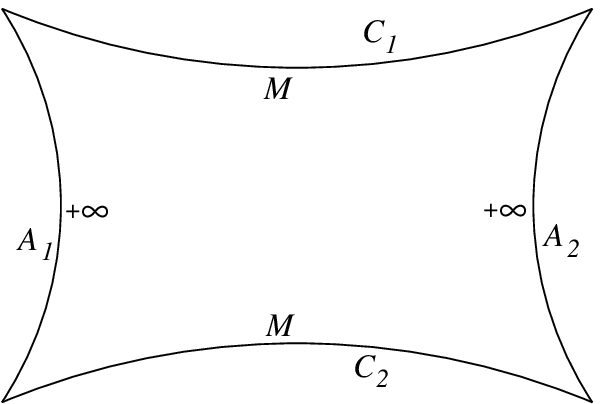}{Figure~3}
\endinsert

If $\vert A_1\vert+\vert A_2\vert< \vert C_1\vert+\vert C_2\vert$, then there is a minimal solution $v$ defined inside ${\cal P}$ taking the boundary values $+\infty$ on $A_1\cup A_2$ and any prescribed continuous data on $C_1\cup C_2$.

Now place the quadrilateral in $\overline D$, putting $C_1$ into $C$, and the rest of ${\cal P}$ inside $D$; cf. figure~4.

\midinsert
\figure{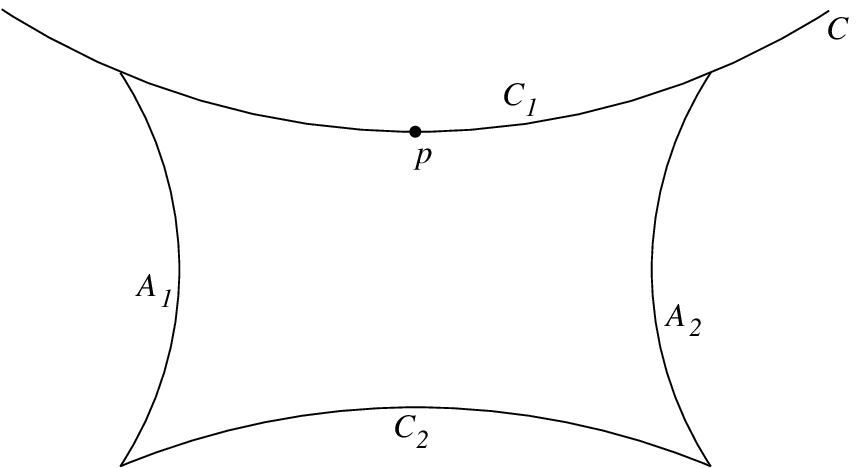}{Figure~4}
\endinsert

Let $M$ be the maximum of $u$ on $C_2$ and $f$ on $C_1$. Let $v$ be the solution in ${\cal P}$ that equals $M$ on $C_1\cup C_2$ and $+\infty$ on $A_1\cup A_2$. Then $v\geq u_n$ on $\partial{\cal P}$ so $v\geq u_n$ in ${\cal P}$. Thus $v$ bounds $u$ in a smaller compact neighborhood of $p$ inside ${\cal P}$.

We complete this section by stating the existence theorem for compact domains.

\theo{Existence Theorem}{Let ${\cal W}$ be a bounded domain with $\partial{\cal W}$ a Jordan curve. Assume there is a finite set $E\subset\partial{\cal W}$ and $\partial{\cal W}-E$ is composed of convex (towards ${\cal W}$) arcs. Then there is a solution to the Dirichlet problem in ${\cal W}$ taking on arbitrarily prescribed continuous data on $\partial{\cal W}-E$. The arcs need not be strictly convex.}

\titre{3.~The Dirichlet Problem on Unbounded Domains}

Let $\Gamma$ be an ideal polygon of $\HM$; i.e., $\Gamma$ is a geodesic polygon all of whose vertices are at infinity $\partial_\infty(\HM)$. We assume $\Gamma$ has an even number of sides $A_1,B_1,A_2,B_2,\ldots,A_k,B_k$, ordered by traversing $\Gamma$ clockwise. Let $D$ be the interior of the convex hull of the vertices of $\Gamma$; so that $\partial D=\Gamma$ and $D$ is a topological disk.

We consider the Dirichlet problem for the minimal surface equation in $D$ where one prescribes the data $+\infty$ on each $A_i$, and $-\infty$ on each $B_j$.

When $D$ is a relatively compact domain in $\RM^2$, with boundary composed of line segments $\{A_i\}$ and $\{B_j\}$, and convex arcs $\{C_l\}$, this type of Dirichlet problem was solved completely by H.~Jenkins and J.~Serrin [\ref{J-S}] (allowing continuous data on the $\{C_l\}$). They found necessary and sufficient conditions -- in terms of the geometry of $\partial D$ -- that guaranty a solution to the Dirichlet problem, for the minimal surface equation in $\RM^2$.

The Jenkins-Serrin theorem was extended to $\HM\times\RM$, for compact domains $\overline D\subset\HM$, and minimal graphs over some non-compact domains were considered [\ref{N-R}].

When $\Gamma$ is an ideal polygon, we will find necessary and sufficient conditions, in terms of the ``lengths'' of the edges, which enable us to solve the Dirichlet problem. Since the lengths of the edges are infinite, we proceed as follows.

At each vertex $a_i$ of $\Gamma$, place a horocycle $H_i$; do this so $H_i\cap H_j=\emptyset$ if $i\neq j$. Let $F_i$ be the convex horodisk with boundary $H_i$.

Each $A_i$ meets exactly two horodisks. Denote by $\tilde A_i$ the compact arc of $A_i$ which is the part of $A_i$ outside the two horodisks, $\partial\tilde A_i$ is two points, each on a horocycle; $\vert A_i\vert$ is the distance between these horocycles, i.e. the length of $\tilde A_i$. Define $\tilde B_i$ and $\vert B_i\vert$ in the same way.

Define
$$\eqalign{
a(\Gamma) &= \sum_{i=1}^k \vert A_i\vert,\hbox{ and }\cr
b(\Gamma) &= \sum_{i=1}^k \vert B_i\vert.\cr}$$
Observe that $a(\Gamma)-b(\Gamma)$ does not depend on the choice of horocycles (assuming they are pairwise disjoint).

\theo{Definition~1}{An ideal geodesic polygon ${\cal P}$ is said to be inscribed in $D$ if the vertices of ${\cal P}$ are among the vertices of $\Gamma$; ${\cal P}$ is a simple closed polygon whose edges are either interior in $D$ or equal to an $A_i$ or $B_j$.}

Notice that the definition of $a(\Gamma)$ and $b(\Gamma)$ extends to inscribed polygons. Also, to each inscribed polygon ${\cal P}$ in $\Gamma$ and a choice of horocycles at the vertices, we associate a Jordan curve $\tilde{\cal P}$ in $\HM$ as follows.

Let $\alpha$ and $\beta$ be geodesic edges of ${\cal P}$ with the same vertex $a_i$. Let $\gamma_i$ denote the compact arc of $H_i$ joining $\alpha\cap H_i$ to $\beta\cap H_i$. Then
$\tilde{\cal P}$ is obtained from ${\cal P}$ by removing the non-compact arcs of $\alpha$ and $\beta$ in $H_i$ and replacing them by $\gamma_i$; cf. figure~5. Denote by $\vert{\cal P}\vert$ the length of the boundary arcs of ${\cal P}$ exterior to the horocycles at the vertices of ${\cal P}$; we call this the truncated length of ${\cal P}$.

\midinsert
\figure{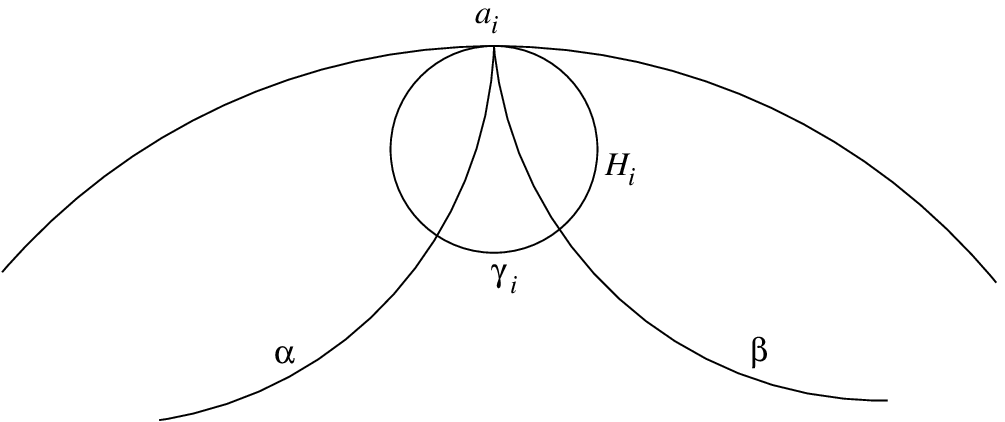}{\rlap{\smash{\raise9mm\hbox{${\cal P}$}}}Figure~5}
\endinsert

We can now state the result.

\theo{Theorem~1}{There is a solution to the minimal surface equation in the polygonal domain $D$, equal to $+\infty$ on each $A_i$ and $-\infty$ on each $B_j$, if and only if the following two conditions are satisfied 
\item{\rm1.}$a(\Gamma)=b(\Gamma)$.\hfill{\tenrm(1)}
\item{\rm2.}for each inscribed polygon ${\cal P}$ in $\Gamma$, ${\cal P}\neq\Gamma$, and for some choice of the horocycles at the vertices, one has
$$2a({\cal P})<\vert{\cal P}\vert\hbox{ and }2b({\cal P})<\vert{\cal P}\vert.\eqno(2)$$
The solution is unique up to an additive constant.}

\rmq{Remark~3}{It is easy to check that if the inequalities of condition~2 hold for some choice of horocycles, then they continue to hold for ``smaller'' horocycles (the inequalities get better); cf. Remark~8. Thus for a given $\Gamma$, a finite number of choices of horocycles suffice to check the inequalities are satisfied.}

\theo{Definition~2}{An inscribed polygon ${\cal P}$ that satisfies the conditions of Theorem~1 is said to be admissible.}

First we fix some notation. At each vertex $a_i$ of $\Gamma$, let $H_i(n)$ be a sequence of nested horocycles at $a_i$, converging to $a_i$ as $n\to\infty$; $F_i(n)$ the convex horodisk bounded by $H_i(n)$.

$H_i(n)\cap\Gamma$ consists of two points and we noted by $\gamma_i(n)$ the compact arc of $H_i(n)$ joining these two points. We now need to work in convex domains so to form the Jordan curves $\hat\Gamma(n)$, we will not attach $\gamma_i(n)$ but the geodesic arc $\hat\gamma_i(n)$, having the same end-points as $\gamma_i(n)$.

Then $\hat\Gamma(n)$ is the convex Jordan curve:
$$\hat\Gamma(n)=[\Gamma-\bigcup_i(\Gamma\cap F_i(n))]\cup[\bigcup_i\hat\gamma_i(n)].$$
Let $D(n)$ denote the disk of $\HM$ bounded by $\hat\Gamma(n)$; so that $\bigcup\limits_n D(n)=D$.

We establish several lemmas for the proof of the Theorem~1.

\theo{Lemma~1}{For each $i=1,\ldots,k$, there is a solution to the Dirichlet problem on $D$ with boundary data $+\infty$ on $A_i$, zero on the rest of $\Gamma$.}

\noindent{\bf Proof:} For each $n$, define $u_n$ to be the minimal graph on $D(n)$, which is equal to $n$ on $\partial D(n)\cap A_i$, and zero on the rest of $\partial D(n)$. By the maximum principle, $u_n$ is a monotone increasing sequence on $D(\ell)$ for $n\geq\ell$. Thus it suffices to show the $u_n$ are uniformly bounded on compact subsets of $D(\ell)-A_i$; then a diagonal process yields the solution to Lemma~1.

In the upper half-plane model of $\HM$ consider the function [\ref{A}], [\ref{SE}]:
$$h(x,y)=\ln\left({\sqrt{x^2+y^2}+y\over x}\right),\ x>0,\ y>0.$$
Consider $h$ defined in the domain of $\HM$ defined by $y>0$, $x>0$; it satisfies the minimal surface equation, has asymptotic values zero on the arc at $\partial_\infty\HM$, given by $x>0$, $y=0$. Also it converges to $+\infty$ on $x=0$, $y>0$.

Thus for any geodesic of $\HM$, one has such a minimal graph. Apply this to $A_i$: let $h$ be defined on the part of $\HM$, with boundary $A_i$, that contains $D$. Clearly $h$ is greater than each $u_n$ by the maximum principle hence on any compact set $K\subset D(\ell)-A_i$, the $u_n$, $n\geq\ell$, are uniformly bounded.

\rmq{Remark~4}{This argument greatly simplifies the proofs of Theorem~2 and Theorem~3, step~1, of [\ref{N-R}]; the barrier $h$ shows the $u_n$ of Theorem~2, and of step~1 of Theorem~3, are uniformly bounded on compact sets.}

\theo{Lemma~2}{For each $i=1,\ldots,k$, there is a solution to the Dirichlet problem on $D$ with boundary data $+\infty$ on each $A_j$, $j\neq i$, and zero on the rest of $\Gamma$.}

\noindent{\bf Proof:} We note by $A=A_1\cup A_2\cup\ldots\cup A_k$, $B=B_1\cup B_2\cup\ldots\cup B_k$, so that $\Gamma=A\cup B$. $\tilde A(n)$ and $\tilde B(n)$ denote the truncated $A$'s and $B$'s in $\hat\Gamma(n)=\partial D(n)$; $\hat\Gamma(n)=\tilde A(n)\cup\tilde B(n)\cup\hat\gamma(n)$, where  $\hat\gamma(n)=\hat\gamma_1(n)\cup\hat\gamma_2(n)\cup\ldots\cup\hat\gamma_k(n)$ is the union of the remaining geodesic arcs. For convenience, suppose $i=1$.

For each $n$, let $u_n$ be the solution to the Dirichlet problem on $D$ with boundary data equal to $n$ on $A_2\cup\ldots\cup A_k$, and zero on the rest of $\Gamma$. This $u_n$ exists: first construct such a solution on each compact $D(\ell)$, then let $\ell\to\infty$ and choose a convergent subsequence by a diagonal process.

By the generalized maximum principle (whose proof we give in Theorem~2), $\{u_n\}$ is a monotone sequence so the Divergence Structure Theorem applies. If the divergence set is not empty, a connected component is an inscribed polygon in $D$ whose geodesic edges join two vertices of $\Gamma$. We show $V=\emptyset$ to prove Lemma~2.

First suppose $V=D$ so $\partial V=\Gamma$.

We fix $\ell$ and consider $n\geq0$, and the function $u_n$ on $D(\ell)$. Let $X_n={\nabla u_n\over W_n}$ where $\nabla u_n$ is the gradient of $u_n$ in $\HM$ and $W_n^2=1+\vert\nabla u_n\vert^2$. Then ${\rm div}(X_n)=0$, so the flux of $X_n$ along $\partial D(\ell)=\hat\Gamma(\ell)$ is zero. For each arc $\alpha\subset\hat\Gamma(\ell)$, let $F_n(\alpha)=\int_\alpha\langle X_n,\nu\rangle\,ds$, where $\nu$ is the outer conormal to $D(\ell)$ along $\partial D(\ell)$.

Then the flux of $X_n$ along $\partial D(\ell)$ yields:
$$0=F_n(\tilde A_1)+F_n(\tilde A_2)+\cdots+F_n(\tilde A_k)+F_n(\tilde B_1)+\cdots+F_n(\tilde B_k)+\sum_{j=1}^k F_n(\hat\gamma_j(\ell)),$$
where $\hat\gamma_j(\ell)$ are the small geodesic arcs in $\partial D(\ell)$. Now the $u_n$ diverge uniformly to infinity on compact subsets of $D$ and are bounded on 
$A_1\cup B_1\cup\ldots\cup B_k$, so
$$\lim_{n\to\infty} F_n(\tilde A_1)=-\vert\tilde A_1\vert,\quad \lim_{n\to\infty} F_n(\tilde B_j)=-\vert\tilde B_j\vert$$
for $j=1,\ldots,k$. Also $F_n(C)\leq\vert C\vert$ for any arc $C\subset\partial D(\ell)$. Hence letting $n\to\infty$, we obtain
$$\vert\hat\gamma(\ell)\vert+\vert\tilde A_2\vert +\cdots+ \vert\tilde A_k\vert\geq\vert\tilde A_1\vert+\vert\tilde B_1\vert+\cdots+\vert\tilde B_k\vert$$
where $\vert\hat\gamma(\ell)\vert$ is the length of all the arcs $\hat\gamma_j(\ell)$ in $\partial D(\ell)$. But for any choice of $\ell$,
$$a(\Gamma)=\vert\tilde A_1\vert +\cdots+ \vert\tilde A_k\vert=\vert\tilde B_1\vert +\cdots+ \vert\tilde B_k\vert =b(\Gamma)$$
by hypothesis and $\vert\hat\gamma(\ell)\vert\to0$, as $\ell\to\infty$ so this last inequality is impossible.

It remains to show $V\neq D$, $V\neq\emptyset$ is impossible. Suppose this were the case. Fix $\ell$ and consider $V\cap D(\ell)=V(\ell)$. $V(\ell)$ is bounded by interior geodesic arcs $\displaystyle C_1,C_2,\ldots$, some arcs $\tilde A_{\displaystyle i_1},\tilde A_{\displaystyle i_2},\ldots,\tilde B_{\displaystyle j_1},\tilde B_{\displaystyle j_2},\ldots$ of $\partial D(\ell)$, and some small geodesic arcs $\displaystyle\tilde\gamma_1,\tilde\gamma_2,\ldots$ (all of these arcs depend on $\ell$).

The flux of $X_n$ along $\partial V(\ell)$ equals zero hence:
$$\sum F_n(\tilde A_i)+\sum F_n(\tilde B_j)+\sum F_n(C_l)+\sum F_n(\tilde\gamma_m)=0$$
where the sums are taken over all the arcs in $\partial V(\ell)$. The arcs $C_1,C_2,\ldots,$ are interior arcs of $D$ so on each $C_l$,
$$\lim_{n\to\infty} F_n(C_l)=-\vert C_l\vert.$$
Similarly, for the $\tilde B_j$ in $\partial V(\ell)$,
$$\lim_{n\to\infty} F_n(\tilde B_j)=-\vert\tilde B_j\vert.$$
Let $\vert\tilde\gamma\vert$ be the total length of the small geodesic arcs in $\partial V(\ell)={\cal P}$, so $\vert\tilde\gamma\vert\geq\sum F_n(\tilde\gamma_m)$. Similarly for the $\tilde A_i$ arcs we have 
$a({\cal P})\geq\sum F_n(\tilde A_i)$.

The flux equality yields
$$a({\cal P})+\vert\tilde\gamma\vert\geq b({\cal P})+\vert C\vert$$
where $\vert C\vert$ is the sum of the lengths of the arcs $C_m$ in ${\cal P}$. This last inequality can be rewritten:
$$2a(P)-a({\cal P})-b({\cal P})-\vert C\vert\geq-\vert\tilde\gamma\vert.$$

Combining the previous inequality and the definition of the truncated length, one obtains
$$2a({\cal P})-\vert{\cal P}\vert\geq-\vert\tilde\gamma\vert.$$
The inscribed polygon ${\cal P}$ satisfies the condition~2 for some choice of horocycles, so $2a({\cal P})-\vert{\cal P}\vert<0$ for this choice of horocycles. For smaller horocycles, this quantity decreases, so the above inequality $2a({\cal P})-\vert{\cal P}\vert\geq-\vert\tilde\gamma\vert$ is impossible for $\ell$ large enough.

\rmq{Remark~5}{Clearly Lemma~2 also proves there is a solution on $D$ which is $-\infty$ on each $B_j$, $j\neq i$, and zero on the rest of the boundary of $D$.}

\noindent{\bf Proof of Theorem~1:} Let $v_n$ be the solution on $D$ which is $n$ on each $A_i$ and zero on each $B_j$. This $v_n$ is easily constructed by solving this Dirichlet problem on the disks $D(\ell)$ and taking the limit as $\ell\to\infty$.

For $0<c<n$, define
$$E_c=\{v_n>c\}\hbox{ and }F_c=\{v_n<c\}.$$

Let $E_c^i$ be the component of $E_c$ whose closure contains $A_i$ and let $F_c^j$ be the component of $F_c$ whose closure contains $B_j$.

We can separate any two of the $A_i$ by a curve joining two of the $B_j$'s, on which $v_n$ is bounded away from $n$. Hence for $c$ sufficiently close to $n$, the sets $E_c^i$ are pairwise disjoint. Let $\mu(n)$ be the infimum of the constants $c$ such that the sets $E_c^i$ are pairwise disjoint.

We claim that each component of $E_c$ is equal to some $E_c^i$. Also, there is a pair $E_\mu^i$ and $E_\mu^j$, $i\neq j$, $\mu=\mu(n)$, such that $E_\mu^i\cap E_\mu^j=\emptyset$ and $\overline{\textstyle\mathstrut E_\mu^i}\cap\overline{\textstyle\mathstrut E_\mu^j}\neq\emptyset$. Hence, given any $F_\mu^i$, there is some $F_\mu^j$ disjoint from it.

To see this claim, consider the level curves $\{v_n=c\}$, for some $0<c<n$. There are no compact level curves by the maximum principle. Each connected level curve is proper hence each end of the level curve is asymptotic to exactly one vertex of $\Gamma$.

Suppose ${\cal W}$ is a component of $E_c$ and ${\cal W}\neq E_c^i$ for each $i$. Then 
${v_n\vrule_{\,_{\scriptstyle{\partial\cal W}}}}=c$ and $\partial{\cal W}$ is composed of level curves of $v_n$ contained in $D$. We know $\partial{\cal W}$ is not compact, so it enters the ends of $D$ and is asymptotic to the vertices at these ends. Form a compact cycle $\beta$ from $\partial{\cal W}$ by attaching short arcs in ${\cal W}$ joining two branches of the curve of $\partial{\cal W}$ asymptotic to the same vertex; cf. figure~6.

\midinsert
\figure{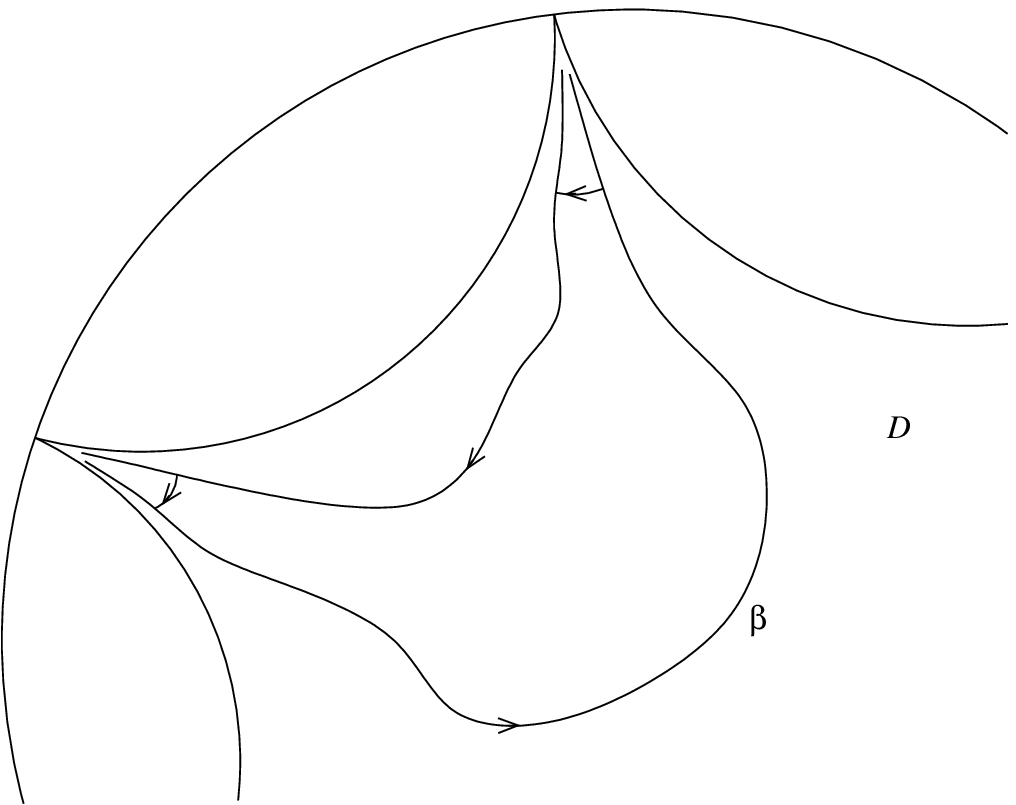}{\rlap{\smash{\raise21mm\hbox{\hskip-14mm$\partial{\cal W}$} 
\raise31mm\hbox{\hskip19mm${\cal W}$}}}Figure~6}
\endinsert

Let $X_n={\nabla v_n\over W_n}$. The flux of $X_n$ along $\beta$ is zero. On the arcs of $\beta$ in $\partial{\cal W}$, the flux is bounded away from zero. Since the flux on the ``short'' arcs of $\beta$ ``inside'' $D$ is bounded by their lengths, and the lengths can be chosen arbitrarily small, this is impossible.

Now is a good time to make a slight digression.

\theo{Theorem~2}{{\rm\bf(Generalized Maximum Principle)} Let $D$ be a domain with $\partial D$ an ideal geodesic polygon. Let ${\cal U}\subset D$ be a domain and $u,v\in C^0(\overline{\cal U})$, two solutions of the minimal surface equation in ${\cal U}$ with $u\leq v$ on $\partial{\cal U}$. Then $u\leq v$ in ${\cal U}$.}

\noindent{\bf Proof:} If this were not true, then we can suppose (after a possible small vertical translation of the graph of $v$) that ${\cal W}=\{p\in {\cal U}\vert u(p)>v(p)\}$ is a non-empty domain with smooth boundary. $\overline{\cal W}$ is not compact by the maximum principle, so ${\cal W}$ has branches going into the ends of $D$. As in the previous argument, we form a compact cycle $\beta$ composed of (long) arcs on $\partial{\cal W}$ and (short) arcs in the ends of $D\cap{\cal W}$; cf. figure~6.

Let $\displaystyle X={\nabla u\over W_u}$ and $\displaystyle Y={\nabla v\over W_v}$. The flux of $X-Y$ across $\beta$ is $0$. On the short arcs of $\beta$, this flux is bounded by twice the length of the short arcs, hence can be made arbitrarily small. It remains to prove the flux of $X-Y$ on the long arcs of $\beta$ is bounded away from zero. This is a well known argument which we repeat here.

Since $u-v>0$ in ${\cal W}$, $u=v$ on $\partial{\cal W}$, we have $\nabla(u-v)=\lambda\eta$ on $\partial{\cal W}$, where $\eta$, the inner pointing conormal to $\partial W$ in $W$, orients the level curve; $\lambda$ is a strictly positive function on $\partial{\cal W}$, by the boundary maximum principle.

An algebraic identity (cf. appendix or [\ref{N-R}]) shows $\langle\nabla(u-v),X-Y\rangle\geq0$ with equality if and only if $\nabla(u-v)$ is zero. Thus $\langle X-Y,\eta\rangle$ is bounded away from zero on the long arcs of $\beta$, which proves the generalized maximum principle.

\rmq{Remark~6}{Theorem~2 extends to possible infinite boundary values for the solutions along geodesic arcs $\alpha$ of $\partial{\cal U}$. It is immediate if the data are distinct by restricting ${\cal U}$. If $u$ and $v$ take the same infinite value along $\alpha$, a careful analysis of $X$ and $Y$ in a neighborhood of $\alpha$ (analogous to the one of the proposition in the appendix) 
gives $X=Y$ along $\alpha$ and the arc $\alpha$ can not generate flux, which allows us to conclude as above. In particular, this proves the unicity part of Theorem~1.}

Now the existence proof of Theorem~1 proceeds as in Jenkins-Serrin [\ref{J-S}].

Let $u_i^+$, $i=1,\ldots,k$ be the solutions which take the value $+\infty$ on $A_i$ and zero on the rest of $\Gamma$. Similarly, let $u_i^-$ be the solution in $D$ which takes the value $-\infty$ on $B_1,\ldots B_{i-1},B_{i+1},\ldots B_k$, and zero on the rest of $\Gamma$; they exist by Lemma~1 and Lemma~2. Define $u_n=v_n-\mu(n)$ and, for each $p\in D$,
$$u^+(p) =\max_i\{ u_i^+(p)\},\ u^-(p)=\min_i\{u_i^-(p)\}.$$

Observe that
$$u^-\leq u_n\leq u^+\hbox{ in }D.$$
To see this, notice that if, for some $p\in D$, $u_n(p)>0$ then $p\in E_\mu^i$ for some $i\in\{1,\ldots,k\}$ and by the generalized maximum principle we have $u_n\leq u_i^+$ in $E_\mu^i$.
If $u_n(p)<0$, then $p\in F_\mu^i$ for some $i$. There is some $j\neq i$ such that $F_\mu^i\cap F_\mu^j=\emptyset$. Then the generalized maximum principle in $F_\mu^i$ yields 
$u_n(p)\geq u_j^-(p)$.

Thus the sequence $u_n$ is uniformly bounded on compact subsets of $D$, so a subsequence of the $u_n$ converges to a solution $u$ in $D$. It remains to show $u=+\infty$ on each $A_i$, and $u=-\infty$ on each $B_j$.

Observe that the sequence of constants $\mu(n)$ diverges to infinity. Otherwise we would have a subsequence (again called $u_n$) such that $u_n$ converges to $u$, $u$ is $+\infty$ on the $A_i$, and $u=-\mu_0$ on the $B_j$.

Fix $\ell$ and calculate the flux of $u$ along $\partial D(\ell)$; this gives
$$\sum_i F_u(\tilde A_i)+\sum_i F_u(\tilde B_i)+\sum_i F_u(\hat\gamma_i(\ell))=0.$$
Here $\tilde A_i$, $\tilde B_i$, $\hat\gamma_i(\ell)$ are the arcs on $\partial D(\ell)=\hat\Gamma(\ell)$; the $\hat\gamma_i(\ell)$ are short geodesic arcs. We have, using the Flux Theorem: 
$$\sum_iF_u(\tilde A_i)=a(\Gamma)=\vert\tilde A_1\vert+\cdots+\vert\tilde A_k\vert,$$
$$\sum_iF_u(\hat\gamma_i(\ell))\leq\vert\hat\gamma(\ell)\vert=\vert\hat\gamma_1(\ell)\vert+\cdots+
\vert\hat\gamma_k(\ell)\vert,$$
and
$$\sum_iF_u(\tilde B_i)\leq b(\Gamma)-\delta\hbox{ for some }\delta>0.$$
Since $\vert\hat\gamma(\ell)\vert\to0$, as $\ell\to\infty$, and $\delta$ can be chosen independent of $\ell$, this contradicts $a(\Gamma)=b(\Gamma)$.

In the same way, one concludes $n-\mu(n)$ diverges to infinity. This completes the proof of Theorem~1.

\titre{4.~Ideal Polygons with Convex Arcs}

Now let $\Gamma$ be an ideal polygon with a finite number of vertices at infinity and geodesic sides $(A_i)_{i=1,\ldots,k}\,$, $(B_j)_{j=1,\ldots,k'}$ as before, but we also allow convex arcs $(C_l)_{l=1,\ldots,k''}$ in $\Gamma$, joining vertices of $\Gamma$ at infinity. We assume the $C_l$ are convex with respect to the domain $D$ bounded by $\Gamma$ and we do not require the $C_l$ to be strictly convex.

We make an important assumption: when a convex arc $C$ in $\Gamma$, has a point $b\in\partial_\infty\HM$ as a vertex, then the other arc $\alpha$ of $\Gamma$ having $b$ as vertex is asymptotic to $C$ at $b$. This means that for a sequence $x_n\in\alpha$, converging to $b$, one has ${\rm dist}_{\HM}(x_n,C)\to 0$ as $n\to\infty$.

This assumption is what we need to assure the Generalized Maximum Principle (Theorem~2) holds in $D$. Its proof is similar to the proof we gave when there are no convex arcs in $\Gamma$.

We now state the expected result.

\theo{Theorem~3}{Let $\Gamma$ be as described above and let $f$ be continuous on the convex arcs $C_l$ of $\Gamma$ ( we assume there are convex arcs $C_l$). Then there is a unique solution in $D$ which is $+\infty$ on each $A_i$, $-\infty$ on each $B_j$ and $f$ on the $C_l$, if and only if
$$2a({\cal P})< \vert{\cal P}\vert,$$
and
$$2b({\cal P})< \vert{\cal P}\vert,$$
for all inscribed polygons ${\cal P}$ in $\Gamma$.}

Recall that we constructed an exhaustion of $D$ by compact convex disks $D(\ell)$ using geodesic arcs at the vertices of $\Gamma$, when there were no convex arcs in $\Gamma$. Now we construct $D(\ell)$ as follows. At each vertex $b$ of $\Gamma$ with an $A_i$ and $C_j$ asymptotic to $b$, choose points $x_\ell\in A_i$ converging to $b$ and let $\gamma(\ell)$ be the minimizing geodesic joining $x_\ell$ to $C_j$. Clearly the $\gamma(\ell)$ are pairwise disjoint. An analogous sequence can be constructed when there are $B$'s and $C$'s at the some vertex $b$ or two arcs $C_i, C_j$ at $b$. Then $D(\ell)$ is defined to have boundary the $\gamma(\ell)$ together with the compact arcs on the $A$'s, $B$'s and $C$'s they bound; cf. figure~7.

\midinsert
\figure{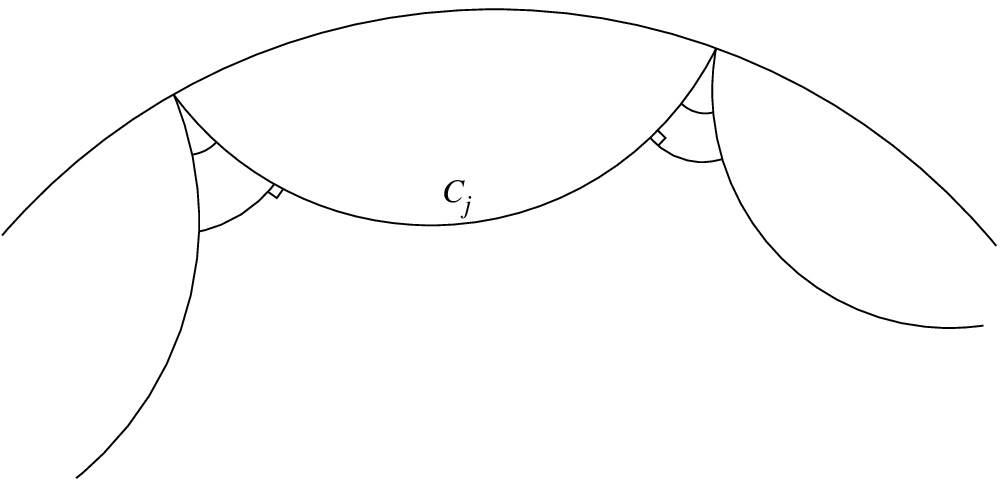}{
\rlap{\smash{
\rlap{\raise37.5mm\hbox{\hskip-30.5mm$x_{\ell'}$}}
\rlap{\raise29.5mm\hbox{\hskip-30.5mm$x_{\ell}$}}
\rlap{\raise29mm\hbox{\hskip-26mm$\gamma(\ell)$}}
\rlap{\raise34.5mm\hbox{\hskip17mm$\gamma(\ell)$}}
\rlap{\raise15mm\hbox{\hskip-10mm$D(\ell)$}}
}}
Figure~7}
\endinsert

First we solve the Dirichlet problem of Theorem~3 when there are no $B$'s, and $f$ is bounded below.

In each $D(\ell)$, let $u_n(\ell)$ be the solution that equals $n$ on the $A_i$ in $\partial D(\ell)$, zero on the $\gamma(\ell)$'s and $\min(n,f)$ on the $C_l$ in $\partial D(\ell)$. By the maximum principle, $u_n(\ell)\leq u_n(\ell')$ on $D(\ell)$ for $\ell\leq\ell'$. Then locally, $u_n(\ell)$ is an increasing bounded sequence as $\ell$ tends to infinity and therefore, by the Monotone Convergence Theorem, it converge to a solution $u_n$ in $D$ that equals $n$ on the $A_i$ and $\min(f,n)$ on the $C_l$.

By the generalized maximum principle, $\{u_n\}$ is a monotone increasing sequence which is uniformly bounded in a neighborhood of strictly convex arcs $C_l$ by the Boundary Values Lemma. Hence if the sequence $\{u_n\}$ has a non-empty divergence set $V$, then a connected component of $V$ is bounded by an inscribed polygon ${\cal P}$ whose edges are geodesics (some of the $C$'s may be in ${\cal P}$) joining vertices of $\Gamma$. Then the same flux calculation we did in Lemma~2 contradicts our hypothesis $2a({\cal P})<\vert{\cal P}\vert$.

Thus the $u_n$ converge to a solution $u$ in $D$ that equals $+\infty$ on the $A_i$ and $f$ on the $C_l$ from the Boundary Values Lemma. The same argument shows that if there are no $A$'s in $\Gamma$ and $f$ is bounded above, then there is a solution in $D$ that equals $-\infty$ on the $B_j$ and $f$ on $C_l$.

Let us now prove the theorem when there are $A$'s, $B$'s and $C$'s.

Let $u^+$ be the solution in $D$ equal to $+\infty$ on the $A_i$, zero on the $B_j$, and $\max(f,0)$ on $C_l$; and similarly, let $u^-$ be the solution in $D$ equal to zero on the $A_i$, $-\infty$ on the $B_j$, and $\min(f,0)$ on $C_l$. These solutions exist by our previous discussion.

Let $v_n$ be the solution in $D$ equal to $n$ on the $A_i$, $-n$ on the $B_j$ and $f_n$ on $C_l$; where $f_n$ is $f$ truncated above by $n$ and below by $-n$.

By the generalized maximum principle,
$$u^-\leq v_n\leq u^+\hbox{ in } D.$$
Therefore the sequence $\{v_n\}$ is uniformly bounded on compact subsets of $D$ so a subsequence converges uniformly on compact sets to a solution $v$ in $D$. By the Boundary Values Lemma $v$ takes on the desired boundary values on $\Gamma$.

\titre{5.~Conformal type}

We present here the main result of this paper concerning harmonic diffeomorphisms. Clearly, any conformal representation of the graphs given by the following theorem gives rise to a diffeomorphism by projection.

\theo {Theorem~3}{In $\HM\times\RM$, there exists entire minimal graphs over $\HM$ which are conformally the complex plane $\CM$.}

\noindent{\bf Proof:} In a first step, we recursively use Proposition~1 and Proposition~2 to construct an exhaustion of $\HM$ by  compact disks $K_n$ with $K_n\subset\inter K_{n+1}$, and a sequence of minimal graphs $u_n$ over $K_n$ (the restriction of ideal Scherk graphs defined over ideal polygonal domains $D_n$ containing $K_n$) satisfying the following:
\item{\it i)}$\parallel u_{n+1}-u_n\parallel_{C^2(K_n)}<\varepsilon_n$, for some sequence $\varepsilon_n>0$, with $\sum_{n=0}^\infty\varepsilon_n<+\infty$,
\item{\it ii)}For each $j$, $0\leq j<n$, the conformal modulus of the annulus in the graph of $u_n$ over the domain $K_{j+1}-\inter K_j$ is greater than one.

For that, let $\varepsilon_n$ be a sequence of positive real numbers such that $\sum_{n=0}^\infty\varepsilon_n<+\infty$. We assume $(D_j,u_j,K_j)$ are constructed for $0\leq j\leq n$ and satisfy the properties {\it i} and {\it ii} we require above.

Using Proposition~2, attach (perturbed) elementary quadrilaterals $E_\tau$, $E'_\tau$, to {\it all\/} of the pairs of sides of $\partial D_n$, to obtain an ideal Scherk graph $u_{n+1}$ over an enlarged polygonal domain $D_{n+1}$. The $E_\tau$, $E'_\tau$ are attached successively to the pairs of sides of $\partial D_n$; the parameter $\tau$ of each pair attached depends on the previous expanded polygons. 

Also, to be sure we will be able to construct an exhaustion, we need to choose carefully the initial non perturbed regular quadrilateral of Proposition~2. For that, consider at each side $\gamma$ of $\partial D_n$, the reflexion across the geodesic orthogonal to $\gamma$ and passing through a fixed point $O\in K_0$, and the special regular quadrilateral along this side $\gamma$ which is invariant by the corresponding reflexion; cf. figure~8. Then use these quadrilaterals to construct the $E_\tau$ and $E'_\tau$.

Moreover, in the Proposition~2, we can choose $\varepsilon$ small enough to get $u_{n+1}$ as close as we want to $u_n$ in the $C^2$-topology on the compact $K_n$, so that properties {\it i} and {\it ii} are satisfied. For we can easily ensure that first 
$\parallel u_{n+1}-u_n\parallel_{C^2(K_n)}<\varepsilon_n$, and secondly, as for $u_n$, the graph of $u_{n+1}$ over each annulus 
$K_{j+1}-\inter K_j$, $0\leq j<n$, has conformal modulus greater than one, since the closer the graphs are, the closer are the conformal moduli.

Now the graph of $u_{n+1}$ is conformally $\CM$ by Proposition~1, so there is a compact $E$ in this graph satisfying:
\item{-}$E$ is a disk that contains $F$ in its interior, and 
\item{-}the conformal modulus of the annulus $E-F$ is greater than one.

\noindent Here $F$ is the graph of $u_{n+1}$ over $\inter K_n$. 

Then define $K_{n+1}$ to be the vertical projection of $E$; eventually enlarge $E$ in order that $K_{n+1}$ has its boundary in a tubular neighborhood of radius one of $\partial D_{n+1}$. By the above construction, this $K_{n+1}$ satisfies property  {\it ii}.
Then the sequence is constructed and the argument will be complete if we prove the $K_n$ exhaust $\HM$.

For each $n$, using the particular geometry of the perturbed quadrilateral we attach to all of the sides of $\partial D_n$, we get that the boundary of $D_{n+1}$ is a fixed constant farther from the fixed point $O$ we have chosen in $K_0$ (if the quadrilateral we add were regular, it would be at least $\ln(1+\sqrt{2})$ by an elementary computation; cf. figure~8, where the equidistant curves $H$ and $H'$ are two parallel horocycles). 
Hence $\partial D_n$ diverge to infinity with $n$. But we constructed $K_n$ so that $\partial K_n$ is uniformly close to $\partial D_n$.

\midinsert
\figure{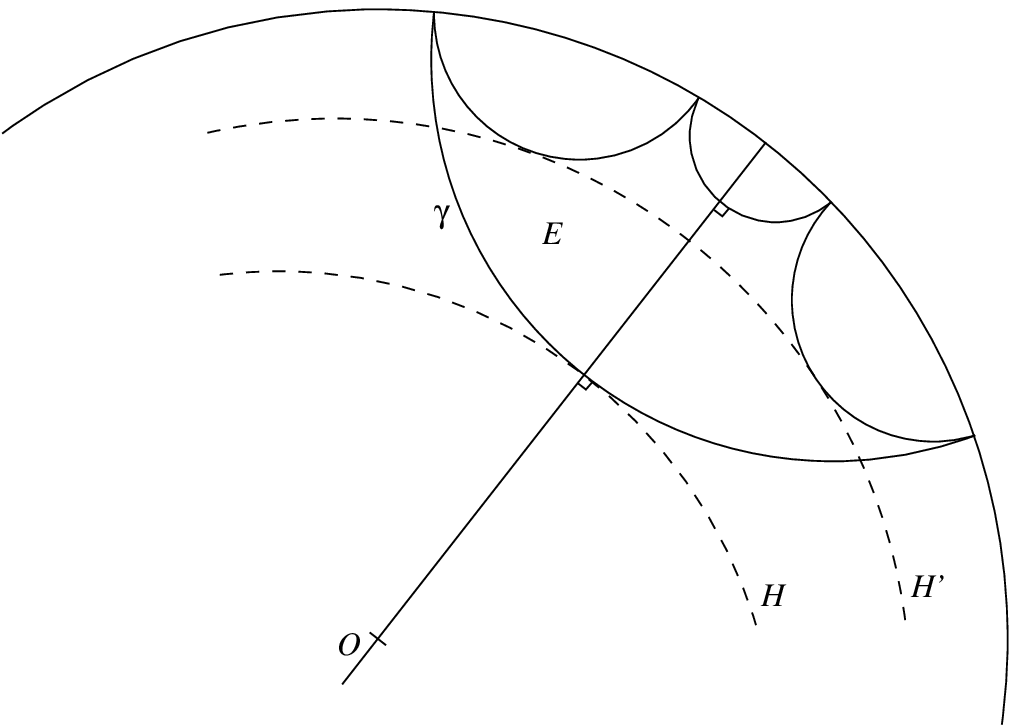}{Figure~8}
\endinsert

Now the second step: we let $n$ tend to infinity. We obtain an entire graph $u$, since at any $x\in\HM$, the $u_n(x)$ form a Cauchy sequence. Since on each $K_{j+1}-\inter K_j$, the $u_n$ converge uniformly to $u$ in the $C^2$-topology, the modulus of the graph of $u$ over $K_{j+1}-K_j$ is at least one. Hence, by Gr\"otzsch Lemma [\ref{V}], the conformal type of the graph of $u$ is $\CM$.

\theo{Proposition~1}{The conformal type of an ideal Scherk graph is $\CM$}.

\noindent{\bf Proof:} Let $\Gamma$ be an ideal Scherk polygon, $D$ the convex hull of $\Gamma$ and let $u$ be a minimal solution defined in $D$ taking the values $+\infty$ on the $A_i$ in $\Gamma$ and $-\infty$ on the $B_j$. Let $\Sigma$ be the graph of $u$.

$\Sigma$ is stable so has uniformly bounded curvature by Schoens' curvature estimates ($\Sigma$ is complete). 

Then there is a $\delta>0$ such that for each $x\in\Sigma$, $\Sigma$ is a graph (in a neighborhood of $x$) over $D_\delta(x)\subset T_x\Sigma$; where $D_\delta(x)$ is the disk of radius $\delta$ in the tangent space $T_x\Sigma$ of $\Sigma$ at $x$, centered at the origin of $T_x\Sigma$. Moreover this local graph has bounded geometry, independent of $x\in\Sigma$. Let $G_\delta(x)$ denote this local graph, $G_\delta(x)\subset\Sigma$.

For $p\in D$, we denote by $\Sigma_\delta(p)$ the local graph $G_\delta(x)$ translated vertically so that $x$ goes to height zero; here $x=(p,u(p))$.

Let $\gamma$ be one of the geodesic components of $\Gamma$ (an $A_i$ or $B_j$) and let $p_n\in D$, $q\in\gamma$ and $\lim\limits_{n\to\infty} p_n=q$.

We claim the local surfaces $\Sigma_\delta(p_n)$ converge uniformly to $\gamma\times\RM$ as $p_n\to q$.

More precisely, they converge uniformly to $\gamma_\delta(q)\times[-\delta,\delta]$, where $\gamma_\delta(q)$ is the $\delta$-interval of $\gamma$ centered at $q$.

Suppose this claim were not the case. First observe the tangent planes to $\Sigma_\delta(p_n)$, at $p_n$, must converge to the vertical plane $\gamma\times\RM$ at $q$. If not, let $q_n$ denote a subsequence of $p_n$ such that the $T_{q_n}\Sigma_\delta(q_n)$ converge to a plane $P$ at $q$, distinct from $Q=\gamma\times\RM$. Since the graphs $\Sigma_\delta(q_n)$
have bounded geometry, for $n$ large, $\Sigma_\delta(q_n)$ is a graph over the $\delta/2$ disk in $P$, centered at $q$; noted $P_{\delta/2}(q)$. A subsequence of these graphs converges uniformly to a minimal graph $F$ over $P_{\delta/2}(q)$. Since $P\neq Q$, there are points of $F$ near 
$q$ whose horizontal projection to $\HM$ is outside of $\overline D$. But the $\Sigma_\delta(q_n)$ converge uniformly to $F$ so for $n$ large, $\Sigma_\delta(q_n)$ would
not be a vertical graph over a domain in $D$; a contradiction.

Thus the tangent planes at $p_n$ to $\Sigma_\delta(p_n)$, converge to the tangent plane to $\gamma\times\RM$ at $q$. A reasoning similar to the previous arguments shows the $\Sigma_\delta(p_n)$ converge to $\gamma_\delta(q)\times[-\delta,\delta]$. A subsequence of the $\Sigma_\delta(p_n)$ converges to a minimal graph $F$ over $\gamma\times\RM$. Were $F$ different
from $\gamma\times\RM$, $F$ would have points near $q$ whose horizontal projection is outside of $\overline D$. But then $\Sigma_\delta(p_n)$ as well, for $n$ large.

Now let $\ell>0$ and suppose $\gamma(\ell)$ is a segment of $\gamma$ of length $\ell$. For $\varepsilon>0$, there exists a height $h=h(\ell,\varepsilon)$ and a tubular neighborhood $T$ of $\gamma(\ell)$ in $\overline D$, such that the graph of $u$ over $T$ is $\varepsilon$-close in the $C^2$-topology, to $\gamma(\ell)\times[h,+\infty)$ when $\gamma$ is an $A_i$, and $\Sigma$ is $\varepsilon$-close over $T$ to $\gamma(\ell)\times(-\infty,h]$ when $\gamma$ is a $B_j$. This follows from our previous discussion by analytic continuation of the disks of radius $\delta$ on $\Sigma$ that converge to $\delta$-disks on $\gamma\times\RM$ as one converges to $\gamma$.

We denote by $\Sigma(\gamma(\ell))$ this part of $\Sigma$, above (or below) height $h$, that is $\varepsilon$-close to $\gamma(\ell)\times [h,\infty)$ (or to $\gamma(\ell)\times(-\infty, h]$). As one goes higher (or lower), the $\Sigma(\gamma(\ell))$ converge to $\gamma(\ell)\times\RM$. In particular, the horizontal projection of $\Sigma(\gamma(\ell))$ to $\gamma(\ell)\times\RM$ is a quasi-isometry.

Now consider a vertex of $\Gamma$ and let $A_i$ and $B_i$ be the edges of $\Gamma$ at this vertex. Let $F_i$ be a horocycle at the vertex and $E_i\subset D$ the inside of the horocycle $H_i$; cf. figure~9.

\midinsert
\figure{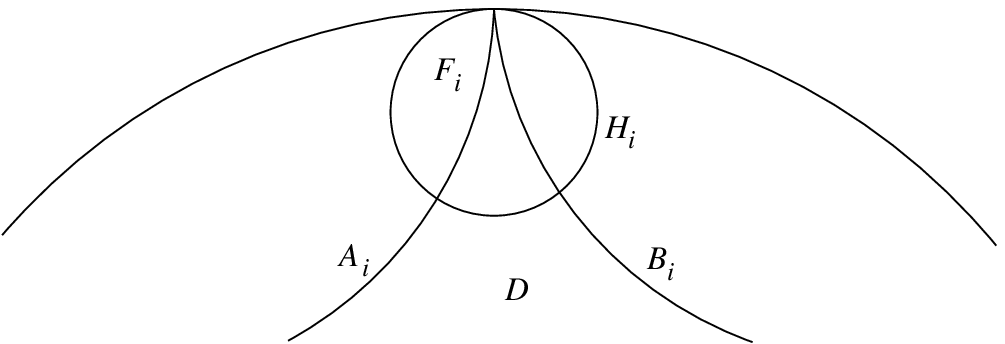}{Figure~9}
\endinsert

By choosing $H_i$ (``small'') we can guarantee that each point of $F_i$ is as close to $A_i\cup B_i$ as we wish. Then for $\varepsilon>0$, there exists an $H_i$ such that the part of $\Sigma$ over $F_i$ (we call $\Sigma(F_i)$) is $\varepsilon$-close to $A_i\times\RM$ (and to $B_i\times\RM$). We choose $\varepsilon$ small so that the horizontal projection of $\Sigma(F_i)$ to $A_i\times\RM$ is a quasi-isometry onto its image.

Let $H_1,\ldots,H_{2k}$ be small horocycles at each of the vertices of $\Gamma$ (we assume $\Gamma$ has $2k$ vertices) so that each $\Sigma(F_i)$ is quasi-isometric to $A_i\times\RM$.

Let $\tilde A_i$ and $\tilde B_i$ denote the compact arcs on each $A_i$ and $B_i$ outside of each $F_i$. For $\vert h\vert$ sufficiently large and $T$ a small tubular neighborhood of $\bigcup\limits_{i=1}^k (\tilde A_i\cup \tilde B_i)$, each component of the part of $\Sigma$ over $T$ projects horizontally to $\tilde A_i\times [h,\infty)$ or to $\tilde B_i\times(-\infty,h]$, quasi-isometrically.

To prove $\Sigma$ is conformally the complex plane $\CM$, we write $\Sigma=\bigcup\limits_{j=0}^\infty K_j$ where each $K_j$ is a disk, $K_j\subset\smash{\inter}K_{j+1}$ for each $j$, and the conformal modulus of each annulus $K_{j+1}-\inter K_j$ is at least one.

Let $K_0$ be the part of $\Sigma$ over $D-\left(T\cup\left(\bigcup\limits_{i=1}^{2k}F_i\right)\right)$. Choose $h_1$ large so that $\Sigma\cap(\HM\times[-h_1,h_1])$ contains an annulus  $K_1-K_0$, $K_1$ compact,  of conformal modulus at least one. Similarly, choose $h_2>h_1$, so that $\Sigma\cap (\HM\times[-h_2,h_2])$ contains an annulus $K_2-K_1$ of conformal modulus at least one. The part of $\Sigma$ outside these $K_j$ is converging to geodesics $\times\RM$, so the $K_j$ exist for all $j$. Thus each ideal Scherk surface is conformally $\CM$.

 \rmq{Remark~7}{We now give another proof that an ideal Scherk graph is conformally $\CM$.  A complete Riemannian surface of finite total curvature is conformally equivalent to a compact Riemann surface punctured in a finite number of points [\ref{Hu}]. Since a Scherk graph is complete and 
simply connected, it suffices to show the graph has finite total curvature.

Recall the notation in the proof of Lemma~2:  
$\hat\Gamma(\ell)=\partial D(\ell)=\tilde A(\ell)\cup \tilde B(\ell)\cup\hat\gamma(\ell)$. For $n$ fixed, we denote by $v_n(\ell)$ the minimal solution on $D(\ell)$, which equals $n$ on $\tilde A(\ell)$ and zero on the rest of the boundary of $D(\ell)$.  As $\ell$ goes to infinity, the $v_n(\ell)$ converge to the function $v_n$ on $D$ which equals $n$ on the $A_i$, $i=1,\ldots,k$, and zero on the $B_j$, $j=1,\ldots.k$. By Gauss-Bonnet formula, the graphs of the $v_n(\ell)$ have total curvature equal to $2\pi-4k{\pi\over2}$. Since the $v_n(\ell)$ converge uniformly on compact sets to $v_n,$ the absolute value of total curvature of the graph of $v_n$ is at most $\vert(2 - 2k)\pi|$.

After renormalization (let $u_n=v_n-\mu_n$), the $u_n$ converge uniformly on compact sets to a Scherk graph $u$ defined on $D$. Thus the graph of $u$ has finite total curvature.  Our analysis of the geometry of the Scherk end shows the total curvature equals $(2 - 2k)\pi$.}

\titre{6.~extending an Ideal Scherk surface}

We now identify $\HM$ with the unit disk in the complex plane 
$\{z\in\CM\ / \ \vert z\vert<1\}$, and $\partial_\infty\HM$ with the unit circle $S^1$. Let 
$(d_1,d_2,\ldots,d_n)$ be n distinct points of $S^1$, ordered clockwise, and denote by 
 $P(d_1,d_2,\ldots,d_n)$ the convex hull of the $n$ points in $\HM$ (for the hyperbolic metric). We say $P(d_1,d_2,d_3,d_4)$ is a regular quadrilateral when the cross ratio 
 $\displaystyle{(d_1-d_3)(d_2-d_4)\over(d_2-d_3)(d_1-d_4)}$ equals $2$. When $n$ is even, we make the convention that ideal Scherk graphs on $P(d_1,d_2,\ldots,d_n)$ take the value $+\infty$ on the geodesic sides $[d_i,d_{i+1}]$, $i$ even, and $-\infty$ on the other sides when $i$ is odd. We denote by ($A_i$) the geodesics $[d_{2i},d_{2i+1}]$, and by $B_i$ the geodesics $[d_{2i+1},d_{2i+2}]$. 

Let ${\cal P}=\partial P(d_1,d_2,\ldots,d_n)$ be an ideal polygon. As in the previous section, place a horocycle $H_i$ at each vertex $d_i$, and let $\displaystyle\vert{\cal P}\vert=\sum_{i=0}^{i=n-1}\vert d_id_{i+1}\vert$ denote the truncated perimeter, where $d_0=d_n$ and $\vert d_id_j\vert$ is the distance between the horocycles $H_i$ and $H_j$. $\vert{\cal P}\vert$ represents the total length of arcs of ${\cal P}$ exterior to all the horocycles. The quantity $\vert{\cal P}\vert$, as the distances $\vert d_id_j\vert$, extends naturally to geodesic polygons with vertices in $\HM$, place horocycles only at vertices which are at infinity. The same extension can be done for the quantities $a({\cal P})$ and $b({\cal P})$ if the polygon ${\cal P}$ comes from a Dirichlet problem.

\rmq{Remark~8}{The utilization of the truncated perimeter $\vert{\cal P}\vert$ gives rise practically to quantities associated to ${\cal P}$ which are independent of the choice of the horocyle at each vertex $d_i$. This allows us to 
check the conditions of Theorem~1 for a choice of horocycles a priori, and a restricted class of inscribed polygons. For, under the hypothesis of Theorem~1, let ${\cal P}$ be an inscribed polygon, ${\cal P}\neq\Gamma$. Notice that, on the one hand at a vertex of ${\cal P}$ with 
a side $A_i$ (necessarily unique), the quantity $\vert{\cal P}\vert-2a({\cal P})$ does not depend on the choice of the horocycle at that point; on the other hand, for the remaining vertices, this quantity increases arbitrarily for a choice of horocycles \og small\fg~enough. Consequently, if the sides $A_i$ alternate on ${\cal P}$, the quantity 
$\vert{\cal P}\vert-2a({\cal P})$ depends only on ${\cal P}$. Thus the polygon ${\cal P}$ will satisfy the first inequality $(2)$ if and only if $\vert{\cal P}\vert-2a({\cal P})>0$. In all the other cases, for a choice of horocycles \og small\fg~enough at vertices where no $A_i$ arrives, the quantity $\vert{\cal P}\vert-2a({\cal P})$ becomes arbitrarily large. The first condition of $(2)$ is automatically satisfied. Similarly, the second inequality is equivalent to checking that $\vert{\cal P}\vert-2b({\cal P})>0$, when the $B_j$ alternate on ${\cal P}$ (then this quantity depends only on ${\cal P}$).}

\theo{Proposition~2}{Let $u$ be an ideal Scherk graph on a polygonal domain 
$D=P(a_1,\ldots,a_{2k})$ and $K$ a compact of $D$. Let 
$D_0=P(b_1,b_2,a_1,b_3,b_4,a_2\ldots,a_{2k})$ be the polygonal domain $D$ to which we attach two regular polygons $E=P(b_1,b_2,a_1,a_0)$ and $E'=P(a_1,b_3,b_4,a_2)$; $E$ is attached to the side $[a_0,a_1]=[a_{2k},a_1]$ of $D$ and $E'$ to the side $[a_1,a_2]$, cf. figure~10.

Then for all $\varepsilon>0$, there exists $(b'_i)_{i=2,3}$ and $v$ an ideal Scherk graph on 
$P(b_1,b'_2,a_1,b'_3,b_4,a_2\ldots,a_{2k})$ such that:
$$\vert b'_i-b_i\vert\leq\varepsilon\hbox{ and }\parallel v-u\parallel_{{\cal C}^2(K)}\leq\varepsilon\ .$$}

\midinsert
\figure{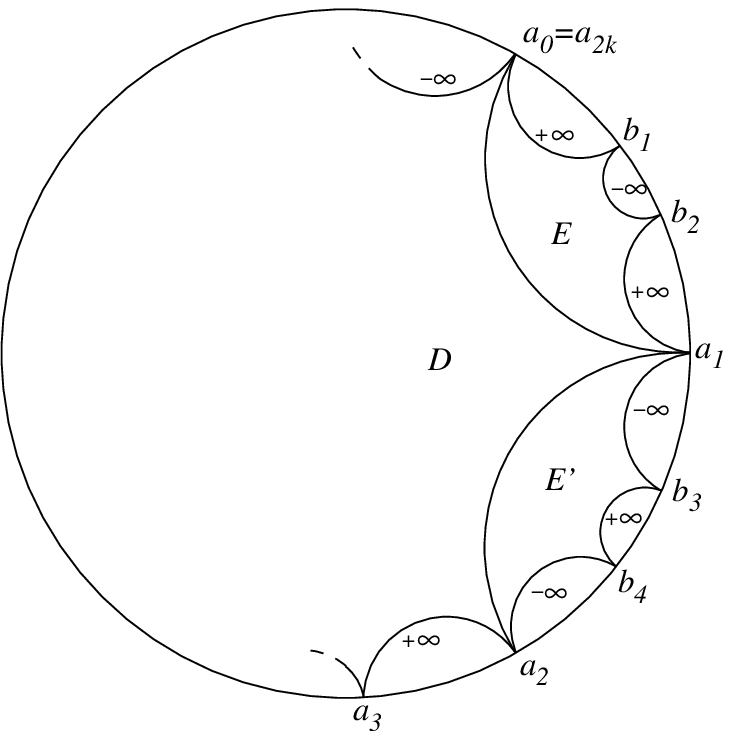}{Figure~10}
\endinsert

We will show that $D_0$ satisfies the conditions of Theorem~1 except for some particular inscribed polygons. This will allow us, by a small variation of $(b_i)_{i=2,3}$, to ensure them completely. Also, it is easy to check that the regularity of $E$ and $E'$ permits the choice, step by step, of the horocycles such that:
$$\eqalign{
&\vert a_0b_1\vert=\vert b_1b_2\vert=\vert b_2a_1\vert=\vert a_1a_0\vert=a,\cr
&\vert a_1b_3\vert=\vert b_3b_4\vert=\vert b_4a_2\vert=\vert a_2a_1\vert=a.\cr}$$
First we establish a Triangle inequality at infinity, using horocycles as usual to define lengths.

\theo{Lemma~3}{{\rm\bf(Triangle inequality at infinity)} For any triangle with vertices $p$, $q$, and $r$ (ideal or not) and small enough pairwise disjoint horocycles placed at the vertices at infinity, we have: 
$\vert pq\vert\leq\vert pr\vert+\vert rq\vert$. If $p$ and $q$ are in $\partial_\infty\HM$, and $r\in\HM$, the inequality is true independently of disjoint horocycles placed at $p$ and $q$.}

\noindent{\bf Proof:} If $r\in\partial_\infty\HM$, the inequality is true for small enough horocycles placed at $r$. If $r$ is in $\HM$, and if $p$ (or $q$) is at infinity, the geodesic $[r,p]$ (or $[r,q]$) is asymptotic to $[p,q]$ at $p$ (or $q$). Then $\forall \alpha>0$, there exists horocycles small enough so that $\vert pq\vert\leq\vert pr\vert+\vert rq\vert+\alpha$. However, the quantity $\vert pr\vert+\vert rq\vert-\vert pq\vert$ does not depend on the 
horocycles placed at $p$ and $q$ if any. Then, passing to the limit, we get a triangle inequality at infinity $0\leq\vert pr\vert+\vert rq\vert-\vert pq\vert$.

In the particular case where $r\in\HM$, $p$ and $q$ at infinity, denote by $H_p$ and $H_q$, the horocycles placed at $p$ and $q$ respectively. If $r$ is in the convex side of one of these horocycles, $H_p$ say, so that $\vert pr\vert=0$, change $H_q$ until touching $H_p$ (on the geodesic $[p,q]$). The Triangle inequality becomes $0\leq\vert rq\vert$ which is true. If $r$ is outside the horocycles, change them, $H_p$ say until $r\in H_p$, and use the previous computation. 

\theo{Lemma~4}{All the inscribed polygons of $D_0$ are admissible except the boundaries of $E$, $E'$ and their complements $D_0\setminus E$, $D_0\setminus E'$.}

\noindent{\bf Proof:} For the entire polygon $D_0$, a simple computation gives that, as for the initial polygon $D$, we have $a(\Gamma_0)=b(\Gamma_0)$, where $\Gamma_0=\partial D_0$. We only prove the inequalities $(2)$ for values 
$+\infty$. Then by symmetry of the problem, we will get the inequalities for the $-\infty$ data on the boundary.

We assume now that ${\cal P}=\partial P$ is an inscribed polygon of $D_0$, 
$P=\partial P(d_1,\ldots,d_n)$ where the $(d_i)_{i=1,\ldots,n}$ are vertices of $D_0$, and moreover that $P\neq D_0$, $P\neq E$, $P\neq E'$, $P\neq D_0\setminus E$, $P\neq D_0\setminus E'$.
By Remark~8, it is enough to prove that we have $\vert{\cal P}\vert-2a({\cal P})>0$ for an a priori choice of disjoint horocycles $H_i$ at $d_i$ ($i=1,\ldots,n$) when the sides $A_i$ where $v=+\infty$ alternate on ${\cal P}$, we assume this hypothesis from now on. 

Moreover, some inequalities used in the proof of Lemma~4 implicitly assume a further choice of horocycles at $d_i$ ($i=1,\ldots,n$). This hypothesis allows us to make that change; cf. Remark~8. 

Consider ${\cal P}'=\partial P$, $P'=P\setminus E'$.

\noindent{\bf Claim: }{\it if $\vert{\cal P'}\vert-2a({\cal P'})>0$ then $\vert{\cal P}\vert-2a({\cal P})>0$.}

\noindent{\bf Proof of the claim:} Consider ${\cal P'}$; if ${\cal P'}={\cal P}$ there is nothing to prove, otherwise the geodesic $[b_3,b_4]$ where $v=+\infty$ is in ${\cal P}$. For convenience, change the notation so that ${\cal P}=\partial P(d_1,b_3,b_4,d_2\ldots,d_n)$. Let $q_1=[d_1,b_3]\cap[a_1,a_2]$ and $q_2=[d_2,b_4]\cap[a_1,a_2]$, cf. figure~11. Notice that if $a_1\in{\cal P}$ (resp. $a_2\in{\cal P}$) then $q_1=a_1$ and
$\vert a_1q_1\vert=0$ by convention (resp. $q_2=a_2$ and $\vert a_2q_2\vert=0$). We have the relations:
$$\eqalign{a({\cal P})&=a({\cal P'})+a,\cr
\vert{\cal P}\vert&=\vert{\cal P'}\vert-\vert q_1q_2\vert+(\vert q_1b_3\vert+a+\vert q_2b_4\vert).\cr}$$
Now, $\vert{\cal P'}\vert-2a({\cal P'})>0$ hence by substitution:
$$\eqalign{\vert{\cal P}\vert-2a({\cal P})&>(\vert q_1b_3\vert+\vert q_2b_4\vert)-(a+\vert q_1q_2\vert)\cr 
&=(\vert q_1b_3\vert+\vert q_2b_4\vert)-(2a-\vert a_1q_1\vert-\vert a_2q_2\vert)\cr
&=(\vert a_1q_1\vert+\vert q_1b_3\vert-a)+(\vert a_2q_2\vert+\vert q_2b_4\vert-a).\cr}$$
 By the Triangle inequality at infinity (directly if $q_1=a_1$ or $q_2=a_2$):
$$\eqalign{\vert a_1q_1\vert+\vert q_1b_3\vert-a&\geq0\cr
\vert a_2q_2\vert+\vert q_2b_4\vert-a&\geq0.\cr}$$
Hence, $\vert{\cal P}\vert-2a({\cal P})>0$ and the claim is proved.

\midinsert
\figure{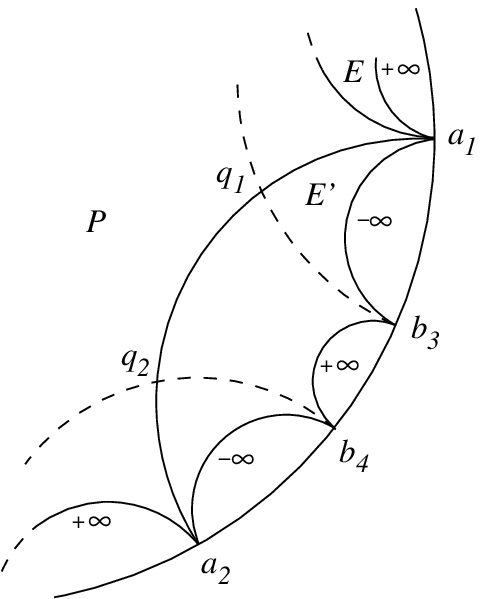}{Figure~11}
\endinsert

So it remains to prove $\vert{\cal P'}\vert-2a({\cal P'})>0$. For that, define 
${\cal P''}=\partial P''$, $P''=P'\setminus E$. 

The key point is a flux inequality for ${\cal P''}$ coming from the initial solution $u$ defined on $D$. We have, $P''\subset D$ and there exists the divergence free field $X={\nabla u\over W}$ 
($W=(1+\vert\nabla u\vert^2)^{1/2}$) on $P''$. Moreover, on the arcs of ${\cal P''}$, $X=\nu$ if $u=+\infty$, $X=-\nu$ if $u=-\infty$, where $\nu$ is the outward normal of $P''$. Let us write ${\cal P''}$ as $I_0\cup I_1 \cup J$: $I_0$ is the union of all geodesics $A_i$ (where $u=+\infty$) contained in ${\cal P''}$ and disjoint from $[a_0,a_1]$, 
$I_1={\cal P''}\cap[a_0,a_1]$ and $J$ the union of the remaining arcs.

The flux of $X$ along ${\cal P''}=\partial P''$ is zero, which yields:
$$0=a({\cal P''})+\vert I_1\vert+F_u(J)+\rho.$$
Here the flux $F_u(J)$ is taken on the compact part of $J$ outside the horocycles and $\rho$ a residual term corresponding to the flux of $X$ along some parts of horocycles. Next we have the truncated perimeter of ${\cal P''}$:
$$\vert{\cal P''}\vert=a({\cal P''})+\vert I_1\vert+\vert J\vert.$$
Adding these last two equalities:
$$\vert{\cal P''}\vert-2(a({\cal P''})+\vert I_1\vert)=\vert J\vert+F_u(J)+\rho.$$
Remark that the condition on the $A_i$ sides of ${\cal P}$ (alternate) yields: 
on the one hand the quantity we have to estimate $\vert{\cal P'}\vert-2a({\cal P'})$ is independent of the choice of horocycles, so this allows us to change this choice as necessary; on the other hand, we have that $P''\neq D$ and $P''\neq \emptyset$. For, if $P''=D$, a careful analysis of the possibilities of inscribed polygons ${\cal P}$ with alternate $A_i$ sides (using $a_0$ and $a_1$ are vertices of ${\cal P}$, cf. figure~10) leads to $P=D_0$ or $P=D_0\setminus E$ which are excluded by hypothesis. Similarly, if $P''=\emptyset$, then $P\subset E$ or $P\subset E'$ and, in this case, the only possibility of an inscribed polygon with alternate $A_i$ sides is $E$ which is excluded too. Therefore $J$ contains interior arcs and as the horocycles at vertices of ${\cal P''}$ diverge:
$$\exists c_0>0\hbox{ so that }\vert J\vert+F_u(J)\geq c_0.$$
We can ensure $\vert\rho\vert<c_0$ for a suitable choice of horocycles. Hence we get the following flux inequality:
$$\vert{\cal P''}\vert-2(a({\cal P''})+\vert I_1\vert)>0.\eqno(3)$$
We have three cases to consider.

\noindent{\bf case~1.} Suppose $[a_0,b_1]\cup[b_2,a_1]\subset{\cal P}$. In this case, $E\subset P'$ and
$$\eqalign{a({\cal P'})&=a({\cal P''})+2a,\cr
\vert{\cal P'}\vert&=\vert{\cal P''}\vert+2a,\cr
\vert I_1\vert&=a.\cr}$$
The flux inequality $(3)$ directly gives:
$$0<(\vert{\cal P'}\vert-2a)-2(a({\cal P'})-a)=\vert{\cal P'}\vert-2a({\cal P'}).$$

\noindent{\bf case~2.} Suppose only one of the $[a_0,b_1]$ or $[b_2,a_1]$ is contained in 
${\cal P}$; $[a_0,b_1]$ say. If we denote $I_1=[a_0,q]$ (cf. figure~12), we have:
$$\eqalign{a({\cal P'})&=a({\cal P''})+a,\cr
\vert{\cal P'}\vert&=\vert{\cal P''}\vert-\vert I_1\vert+a+\vert b_1q\vert.\cr}$$
The flux inequality $(3)$ yields:
$$\eqalign{0&<(\vert{\cal P'}\vert+\vert I_1\vert-a-\vert b_1q\vert)-2(a({\cal P'})-a+\vert I_1\vert),\cr
0&<\vert{\cal P'}\vert-2a({\cal P'})-\vert b_1q\vert-\vert I_1\vert+a}.$$
Hence, using the Triangle inequality at infinity (also valid if $q=a_1$, which does not occur), we obtain:
$$\vert{\cal P'}\vert-2a({\cal P'})>\vert b_1q\vert+\vert qa_0\vert-a\geq0.$$

\midinsert
\figure{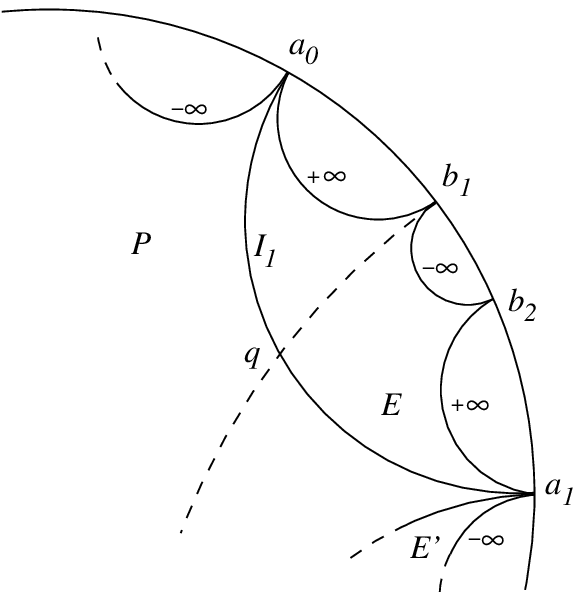}{Figure~12}
\endinsert

\noindent{\bf case~3.} The remaining case is for $P'\subset D$. Then the flux inequality $(3)$ gives directly the result for ${\cal P'}=\partial P'$. This completes the proof of Lemma~4.

The only obstructions to the existence of an ideal Scherk graph on $D_0$ come from the polygons $E$, $E'$ and their complements where we have some cases of equality in $(2)$. In the next lemma, we will ensure them by a perturbation of $D_0$.

\theo{Lemma~5}{There exists $\tau_0>0$ such that for all $\tau\in(0,\tau_0]$, there exist $v_\tau$, an ideal Scherk graph on $D_\tau=P(b_1,b_2(\tau),a_1,b_3(\tau),b_4,a_2,\ldots,a_{2k})$ with:
$\vert b_i(\tau)-b_i\vert\leq\tau$ for $i=2,3$.}

\noindent{\bf Proof:} Let $\tau>0$ and consider moving $b_2$ and $b_3$ in the direction of $a_1$, such that $\vert b_i(\tau)-b_i\vert\leq\tau$ ($i=2,3)$, cf. figure~13. By such variations, the quantities 
$(\vert a_0a_1\vert-\vert a_1b_2(\tau)\vert+\vert b_2(\tau)b_1\vert-\vert b_1a_0\vert)$ and 
$(\vert a_2a_1\vert-\vert a_1b_3(\tau)\vert+\vert b_3(\tau)b_4\vert-\vert b_4a_2\vert)$, which are 
independent of the choice of horocycles, increase. They are zero for the initial polygons $E$ et $E'$. Then there exists such variations $b_2(\tau)$ and $b_3(\tau)$ so that, for $\tau>0$:
$$\eqalign{\varphi(\tau)&=\vert a_0a_1\vert-\vert a_1b_2(\tau)\vert+\vert b_2(\tau)b_1\vert-\vert b_1a_0\vert>0,\cr\varphi(\tau)&=\vert a_2a_1\vert-\vert a_1b_3(\tau)\vert+\vert b_3(\tau)b_4\vert-\vert b_4a_2\vert>0.\cr}$$

\midinsert
\figure{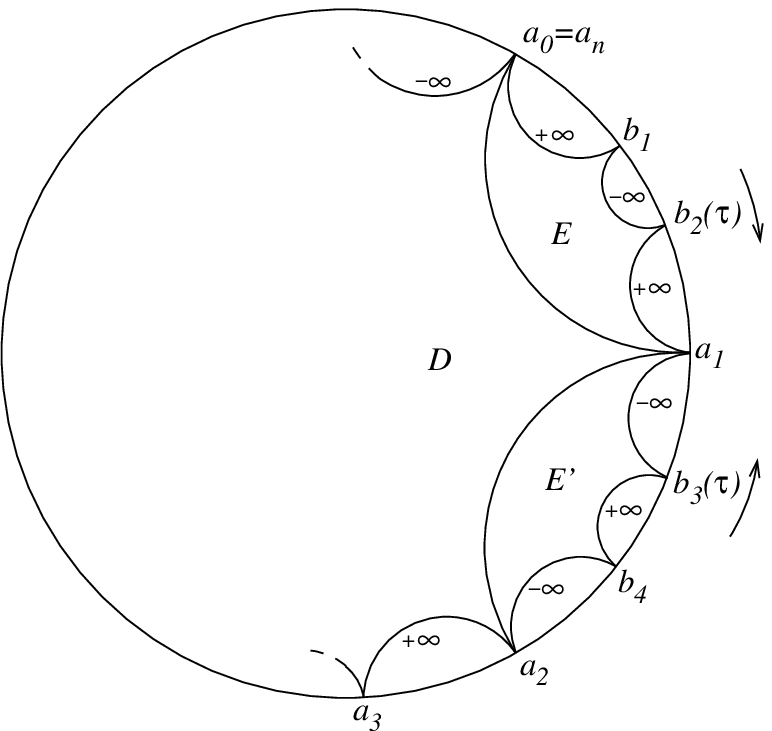}{Figure~13}
\endinsert

Hence, for the polygon $\Gamma_\tau=\partial D_\tau$, we have the condition (1) of Theorem~1, $a(\Gamma_\tau)=b(\Gamma_\tau)$ for all $\tau$ small. On the other hand, for this choice of variation, the $\varphi(\tau)$-perturbed polygons $E_\tau=P(b_1,b_2(\tau),a_1,a_0)$ and $E'_\tau=P(a_1,b_3(\tau),b_4,a_2)$ satisfy
$$\vert\partial E_\tau\vert-2a(\partial E_\tau)=\vert\partial E'_\tau\vert-2b(\partial E'_\tau)=\varphi(\tau)>0,$$
and the inequalities $(2)$ are true for these polygons (use Remark~8). A similar computation gives the inequalities $(2)$ for $D_\tau\setminus E_\tau$ and $D_\tau\setminus E'_\tau$.

In order to prove the inequalities $(2)$ for all other inscribed polygons of $D_\tau$, we use the Lemma~4. For the inscribed polygons ${\cal P}$ of $D_0$ (except $\partial D_0$ and those excluded by Lemma~4), the inequalities $(2)$ are strict and so, are stable by small enough perturbations of vertices (and attached horocycles). There are a finite number of such admissible polygons; thus there exists $\tau_0>0$ such that for all $0<\tau\leq\tau_0$, and variations $b_2(\tau)$ and $b_3(\tau)$ as above, the conditions of Theorem~1 are satisfied. This ensures the existence of $v_\tau$ on $D_\tau$.

\noindent{\bf Proof of Proposition~2:} The main step of the proof consists in establishing that 
$\displaystyle\lim_{\tau\to0}{\nabla v_{\tau}\vrule_{\,_{\scriptstyle D}}}=\nabla u$. For that, consider the divergence free fields associated to $v_\tau$ and $u$: $\displaystyle X_\tau={\nabla v_\tau\over W_\tau}$ and 
$\displaystyle X={\nabla u\over W}$ (with $W_\tau=(1+\vert\nabla v_\tau\vert^2)^{1/2}$, 
$W=(1+\vert\nabla u\vert^2)^{1/2}$); we will prove that 
$\displaystyle\lim_{\tau\to0} {X_{\tau}\vrule_{\,_{\scriptstyle D}}}=X$.

On the boundary, let $\nu$ be the outer pointing normal to $\partial D$. On 
$\partial D\setminus([a_0,a_1]\cup[a_1,a_2])$, $v_\tau$ and $u$ take the same infinite values. Hence 
$X_\tau=X=\pm\nu$. On $[a_0,a_1]$, $u=+\infty$ so $X=\nu$. On the other side, consider the boundary of the domain $E_\tau$ truncated by horocycles. Denote the four horocycle arcs by $\tilde\gamma$. An estimate of the flux of $X_\tau$ yields:
$$0=\vert a_0b_1\vert-\vert b_1b_2(\tau)\vert+\vert b_2(\tau)a_1\vert
+\int_{[a'_0,a'_1]}\langle X_\tau,(-\nu)\rangle\,ds+F_{v_\tau}(\tilde\gamma),$$
the integral is on $[a'_0,a'_1]$, the compact part of $[a_0,a_1]$ joining the horocycles. Then 
$$0=-\varphi(\tau)+\int_{[a'_0,a'_1]}(1-\langle X_\tau,\nu\rangle)\,ds+F_{v_\tau}(\tilde\gamma).$$ 
For a diverging sequence of nested horocycles, we get the convergence of the integral on the whole geodesic and the equality:
$$\int_{[a_0,a_1]}(1-\langle X_\tau,\nu\rangle)\,ds=\varphi(\tau).$$
In the same way, on $[a_1,a_2]$ we get a convergent integral
$$\int_{[a_1,a_2]}(1+\langle X_\tau,\nu\rangle)\,ds=\varphi(\tau).$$
Because of the value of $X$ on $[a_0,a_1]\cup[a_1,a_2]$, for any family $\alpha$ of disjoint arcs of $\partial D$
$$\left\vert\int_\alpha\langle X-X_\tau,\nu\rangle\,ds\right\vert
\leq\int_{[a_0,a_1]\cup[a_1,a_2]}\left\vert\langle X-X_\tau,\nu\rangle\right\vert\,ds
=2\varphi(\tau).
\eqno(4)$$ 

For the study of the field $X_\tau-X$ on the interior of $D$, we consider the flux of $X_\tau-X$ along a level curve through an interior point $p$. This level curve goes to the boundary of $D$ where we create a closed cycle by attaching short curves and a curve on the boundary of $D$ to the level curve. Then the flux is zero along the closed cycle, and small along the curve we attach to the level curve. Thus the flux is small along the level curve, which implies the tangent planes are close, then bounded curvature of the graphs gives the solutions are close, we now make this precise.

As in the previous section, $\Sigma$ the graph of $u$ and $\Sigma_\tau$ the graph of $v_\tau$ are stable, complete and satisfy uniform curvature estimates. Then 
$$\displaylines{\forall\mu>0,\ \exists\rho>0\hbox{ (independent of $\tau$) such that }\forall p\in D:\cr
q\in\Sigma_\tau\cap B((p,v_\tau(p)),\rho)\Longrightarrow
\parallel n_\tau(q)-n_\tau(p)\parallel\leq\mu.\cr}$$
Here, $n_\tau$ denotes the normal to $\Sigma_\tau$ pointing down and $B((p,v_\tau(p)),\rho)$ the ball of radius $\rho$, centered at $(p,v_\tau(p))\in\HM\times\RM$. We have the same estimates for $\Sigma$.

Fix any $\mu>0$ and $p\in D$, this gives a $\rho_1\leq\rho/2$ (independent of $\tau$) such that 
$\forall q\in D(p,\rho_1)$, the disk of $\HM$ with center p and radius $\rho_1$, we have 
$\vert u(q)-u(p)\vert\leq\rho/2$.

Assume now that $\parallel n(p)-n_\tau(p)\parallel\geq3\mu$. Let $\Omega_\tau(p)$ be the connected component of $\{u-v_\tau>u(p)-v_\tau(p)\}$ with $p$ in its boundary, and $\Lambda_\tau$ the component of $\partial\Omega_\tau$ containing $p$. $\Lambda_\tau$, as level curve of $u-v_\tau$, is piecewise smooth. Above $\Lambda_\tau\cap D(p,\rho_1)$, there are two parallel curves: 
$\sigma\subset\Sigma$ and $\sigma_\tau\subset\Sigma_\tau$. Moreover on $\sigma$:
$$\forall q\in\Lambda_\tau\cap D(p,\rho_1),\ \vert(q,u(q))-(p,u(p))\vert\leq\rho_1+\rho/2\leq\rho$$
Hence $\displaystyle\parallel n(q)-n(p)\parallel\leq\mu$.

By a vertical translation of height 
$(v_\tau(p)-u(p))$:
$$(q,v_\tau(q))\in\sigma_\tau\hbox{ and }v_\tau(q)-v_\tau(p)=u(q)-u(p).$$
Then
$$\parallel(q,v_\tau(q))-(p,v_\tau(p))\parallel\leq\rho\hbox{ and }
\parallel n_\tau(q)-n_\tau(p)\parallel\leq\mu.$$
Combining the two last estimates with the assumption on the normals at $p$, we obtain:
$$\forall q\in\Lambda_\tau\cap D(p,\rho_1),\ \parallel n(q)-n_\tau(q)\parallel\geq\parallel n(p)-n_\tau(p)\parallel-2\mu\geq\mu.$$
Apply Lemma~A.1 bellow to conclude 
$$\int_{\Lambda_\tau\cap D(p,\rho_1)}\langle X-X_\tau,\eta\rangle\,ds\geq{\rho_1\mu^2\over2}.$$
As in the proof of the maximum principle, $\langle X-X_\tau,\eta\rangle$ is non negative outside the isolated points where 
$\nabla(u-v_\tau)=0$; then for all compact arcs $\beta\subset\Lambda_\tau$, containing 
$\Lambda_\tau\cap D(p,\rho_1)$ we have:
$$\int_\beta\langle X-X_\tau,\eta\rangle\,ds\geq{\rho_1\mu^2\over2}.\eqno(5)$$

By the maximum principle, $\Lambda_\tau$ is non compact in $D$; so its two infinite branches go close to $\partial D$. Then there exists a connected compact part $\beta$ of $\Lambda_\tau$, containing 
$\Lambda_\tau\cap D(p,\rho_1)$, and two arcs $\tilde\gamma$ in $D$ small enough and joining the extremities of $\beta$ to $\partial D$. Eventually truncating by a family of horocycles, the flux formula for $X-X_\tau$ yields:
$$0=\int_\beta\langle X-X_\tau,(-\eta)\rangle\,ds+\int_\alpha\langle X-X_\tau,\nu\rangle\,ds+
F_{u-v_\tau}(\tilde\gamma\cup\tilde\gamma'),$$
where $\alpha$ is contained in $\partial D$ and $\tilde\gamma'$ is contained in the horocycles and correctly oriented. Using $(4)$ and $(5)$ we obtain
$${\rho_1\mu^2\over2}\leq2\varphi(\tau)+F_{u-v_\tau}(\tilde\gamma\cup\tilde\gamma').$$
When the length of $\tilde\gamma\cup\tilde\gamma'$ goes to zero, we conclude
$${\rho_1\mu^2\over2}\leq2\varphi(\tau).$$
Hence, 
$$\displaystyle\varphi(\tau)\leq{\rho_1\mu^2\over4}\Longrightarrow
\parallel X(p)-X_\tau(p)\parallel\leq\parallel n(p)-n_\tau(p)\parallel\leq3\mu.$$
This gives precisely the behavior of $X_\tau$ and $v_\tau$ for $\tau$ close to zero. After the renormalisation $v_\tau(p_0)=u(p_0)$ for a fixed $p_0\in D$, we have 
$\displaystyle\lim_{\tau\to0} {v_{\tau}\vrule_{\,_{\scriptstyle D}}}=u$. 
The convergence being uniform and ${\cal C}^\infty$ on compact sets, for $\tau$ small enough we can ensure $\displaystyle\parallel v_\tau-u\parallel_{{\cal C}^2(K)}\leq\varepsilon$.
Using Lemma~5, the existence of the ideal Scherk graph of Proposition~2 is established.

We present here the geometric estimate of flux along level curves, we need at the end of the proof of Proposition~2.

\theo{Lemma~A.1}{Let $w$ and $w'$ be two minimal graphs of $\HM\times\RM$, above a domain $\Omega\subset\HM$; $n$ and $n'$ their respective normals. Then at any regular point of $w'-w$:
$$\langle X'-X,\eta\rangle_{\HM}\geq{\parallel n'-n\parallel^2\over4}\geq{\vert X'-X\vert^2\over4},$$
where $X$ (resp. $X'$) is the projection of $n$ (resp. n') on $\HM$ and 
$\displaystyle\eta={\nabla(w'-w)\over\vert\nabla(w'-w)\vert}$ orients the level curve at this regular point.}

\noindent{\bf Proof:} We write $X={\nabla w\over W}$, $X'={\nabla w'\over W'}$ (the normals point down). Classically [\hbox{\ref{Co-K}}],
$$\eqalign{\langle X'-X,\nabla w'-\nabla w\rangle_{\HM}&=\langle n'-n,W'n'-Wn\rangle_{\HM\times\RM}\cr
&=(W+W')\big(1-\langle n,n'\rangle_{\HM\times\RM}\big)\cr
&=(W+W'){\parallel n'-n\parallel^2\over 2}.\cr}$$
Also $\displaystyle{\vert\nabla w'-\nabla w\vert\over W'+W}
\leq{\vert\nabla w'\vert\over W'}+{\vert\nabla w\vert\over W}\leq2$.

\noindent Hence 
$\displaystyle\langle X'-X,\eta\rangle={W+W'\over\vert\nabla(w'-w)\vert}\ {\parallel n'-n\parallel^2\over 2}
\geq{\parallel n'-n\parallel^2\over4}$.

\noindent The last inequality of Lemma~A.1 simply arises by projection.

\titre{References}
\item{[\ref{A}]}U.~Abresch: Personal Communication.
\item{[\ref{Co-K}]}P.~Collin, R.~Krust: {\it Le probl\`eme de Dirichlet pour l'\'equation des surfaces minimales sur des domaines non born\'es}; Bull. Soc. Math. France {\bf119} (1991), p.443--462.
\item{[\ref{He}]}E.~Heinz: {\it \"Uber die L\"osungen der Minimalfl\"achengleichung}; Nachr. Akad. Wiss. G\"ottingen. Math.-Phys. Kl. {\bf1952} (1952), p.51--56.
\item{[\ref{Hu}]}A.~Huber: {\it On subharmonic functions and differential geometry in the large}; Comment. Math. Helv. {\bf 32} (1957), p.13--72.
\item{[\ref{J-S}]}H.~Jenkins, J.~Serrin: {\it Variational problems of minimal surface type. II. Boundary value problems for the minimal surface equation}; Arch. Rational Mech. Anal. {\bf21} (1966), p.321--342.
\item{[\ref{M}]}V.~Markovic: {\it Harmonic diffeomorphisms and conformal distortion of Riemann surfaces};  Comm. Anal. Geom. {\bf10} (2002), p.847--876.
\item{[\ref{N-R}]}B.~Nelli, H.~Rosenberg: {\it Minimal surfaces in $\HM^2\times\RM$}; Bull. Braz. Math. Soc. (N.S.) {\bf33} (2002), p.263--292.
\item{[\ref{SE}]}R.~Sa~Earp: {\it Parabolic and hyperbolic screw motion surfaces in $\HM^2\times\RM$}; to appear in J. Aust. Math. Soc.
\item{[\ref{S}]}R.~Schoen: {\it The role of harmonic mappings in rigidity and deformation problems};Complex Geometry Proceedings of the Osaka Intern. Conf. 1988,  edited by Gen Komatsu and Yusuke Sakane, Marcel Decker Inc. Lecture Notes in Pure and Applied Mathematics {\bf143}, p. 179--200.
\item{[\ref{S-Y}]}R.~Schoen, S.T.~Yau: {\it Lectures on harmonic maps}; Conference Proceedings and Lecture Notes in Geometry and Topology, International Press, Cambridge MA. 1997.
\item{[\ref{V}]}A.~Vassil'ev: {\it Moduli of families of curves for conformal and quasiconformal mappings}; Lecture Notes in Mathematics {\bf 1788} (2002), Springer-Verlag, Berlin.
\bigskip
\begingroup\obeylines\parskip=0mm
{\pcap Pascal Collin:}
Institut de Math\'ematiques de Toulouse, France
collin@math.ups-tlse.fr
\medskip
{\pcap Harold Rosenberg:}
Institut de Math\'ematiques de Jussieu, France
rosen@math.jussieu.fr
\endgroup
\end